\documentclass[12pt]{article}
\usepackage{amsmath}
\usepackage{graphicx,psfrag,epsf}
\usepackage{enumerate}
\usepackage{natbib}
\usepackage{url} % not crucial - just used below for the URL

\newcommand{\blind}{1}
\usepackage{appendix}
\usepackage{amsmath,amsthm}
\usepackage{amsfonts}
\usepackage{epsfig}
\usepackage{graphicx}
\usepackage{epsfig, color, graphicx}
\usepackage{natbib}
\usepackage{mathrsfs}
\usepackage{bm}
\usepackage{multirow}
\usepackage{tikz}
\usepackage{subfigure}
\usepackage{amssymb}%??????
\usepackage{mathrsfs}%????

\newcommand{\bbB}{{\bf B}}

\newcommand{\bbx}{{\bf x}}

\newcommand{\bbX}{{\bf X}}

\newcommand{\bbz}{{\bf z}}
\newcommand{\bbS}{{\bf S}}

\newcommand{\bbI}{{\bf I}}

\newcommand{\bbF}{{\bf F}}

\newcommand{\bqn}{\begin{eqnarray*}}
	\newcommand{\eqn}{\end{eqnarray*}}

\newcommand{\bbA}{{\bf A}}

\newcommand{\bqa}{\begin{eqnarray}}
\newcommand{\eqa}{\end{eqnarray}}

\newcommand{\al}{\alpha}
\newcommand{\la}{\lambda}
\newcommand{\mL}{{\mathscr L}}
\newcommand{\E}{{\mathbb E}}
\renewcommand{\(}{\left(}
\renewcommand{\)}{\right)}

\numberwithin{equation}{section}
\newtheorem{theorem}{Theorem}[section]

\newtheorem{remark}[theorem]{Remark}
\newtheorem{corollary}[theorem]{Corollary}
\newtheorem{lemma}[theorem]{Lemma}

\begin{document}

\def\spacingset#1{\renewcommand{\baselinestretch}%
{#1}\small\normalsize} \spacingset{1}

%%%%%%%%%%%%%%%%%%%%%%%%%%%%%%%%%%%%%%%%%%%%%%%%%%%%%%%%%%%%%%%%%%%%%%%%%%%%%%

\if1\blind
{
  \title{\bf Optimal modification of the LRT for the equality of two high-dimensional covariance matrices}
  \author{Qiuyan Zhang, Jiang Hu\footnote{Corresponding author.}~ and  Zhidong Bai\\ \hspace{.2cm}\\
    KLASMOE \& 
     School of Mathematics and Statistics\\
      Northeast Normal University,  China.}
  \maketitle
} \fi

\if0\blind
{
  \bigskip
  \bigskip
  \bigskip
  \begin{center}
    {\LARGE\bf Optimal modification of the LRT for the equality of two high-dimensional covariance matrices}
\end{center}
  \medskip
} \fi

\begin{abstract}
{This paper considers  the  optimal modification of the likelihood ratio test (LRT) for the equality of two high-dimensional covariance matrices. The optimality  here means that  the modification of LRT cannot be improved anymore in our model settings.  It is well-known that the classical  log-LRT is not well defined when the dimension is larger than  or equal to the sample size. Or even the log-LRT is well-defined, it is usually perceived as a bad statistic in high dimension cases for their low powers under some alternatives. 
 In this paper, we shall argue some goodnesses of the modified log-LRT, and  an optimally modified test that works well in cases where the dimension is  larger than the sample sizes is proposed. Besides, the test is established under the weakest conditions on the moments and the dimensions of the samples. The asymptotic  distribution of the proposed test statistic is also obtained  under  the null  hypotheses. What is more, we also propose  a lite  version of  the modified LRT in the paper.  A simulation study and a real data analysis  show that the  performances of  the two proposed statistics are confirmed to be  invariant to affine transformations. }
\end{abstract}

\noindent%
\vfill

\newpage
\spacingset{1.45} % DON'T change the spacing!

\section{Introduction}

Since  the assumption of homogeneity of covariance matrices is needed in  many multivariate statistical analyses based on  two populations,
the equality of two   covariance matrices is among the most  active hypothesis tests. These tests
date back to the work of  \cite{Wilks32C}, which was followed by a huge amount of literature.  Suppose we have $N:=N_1+N_2$  observations $\{\bbz_i^{(l)}\sim {N}(\bm\mu_l, \Sigma_l), i=1,\dots N_l, l=1,2\}$ and wish  to test the hypothesis
\begin{align}\label{hypo}
H_0:~\Sigma_1=\Sigma_2\quad\mbox{v.s.}\quad H_1:~\Sigma_1\neq \Sigma_2,
\end{align}
where $\bm\mu_l$ and $\Sigma_l$ are the population mean vectors and covariance matrices of the $p$-dimensional vectors $\bbz_i^{(l)}$, $l=1,2$ respectively. It is natural  to
first consider the  likelihood ratio test (LRT) if it is ``applicable". But when is the LRT applicable to testing the equality of two   covariance matrices?  The traditional viewpoint is that the LRT is applicable  if the sample sizes are both much larger than the  dimensions based on the $\chi^2$ approximation  \citep{Wilks46S}.  However, \cite{BaiJ09C} and  \cite{JiangY13C} showed that  when the dimensions are large but smaller than the sample sizes,   the traditional $\chi^2$ approximation of the LRT fails to work well. This problem encouraged us to investigate the   conditions under which the high-dimensional LRT  is applicable. 
%{\color{blue}Therefore, the most important contribution, in this paper, to an applicable high-dimensional LRT is that we present some  so-called ``optimal" conditions. ``Optimal" conditions contain two aspects: the condition on dimensionality which is exactly the ratio of the dimensions to the sample sizes and the moments condition. There will be a specific exposition in Section 2.}
% of  high dimensional LRT for equality of two   covariance matrices,
Additionally,  we   consider a lite LRT  proposed by  \cite{PillaiJ67P,PillaiJ68P}.

Currently, there are  three general types of test procedures for the high-dimensional hypothesis $\eqref{hypo}$  that are widely discussed in the literature: (i) Corrected classical LRTs, see, e.g., \cite{BaiJ09C,JiangJ12L,JiangY13C}; (ii) Nonparametric methods, see, e.g., \cite{LedoitW02S,LiC12T,SrivastavaY10T}; (iii) Maximum element methods, see, e.g., \cite{CaiL13T,CaiM13O}.  The strengths and weaknesses of these three methods are significant.  Nonparametric methods and maximum element methods  can address cases where the dimensions are much larger  than the sample sizes  but  are strongly restricted by the structure or eigenvalues of the population covariance matrices. By contrast,
corrected classical LRT requires  the dimensions to be smaller than the  sample sizes, but  there is no assumption on the population covariance matrices.
In addition,  if the dimensions  are fixed,  LRT has been shown  to be  unbiased and uniformly most
powerful among affine-transform-invariant tests.
Therefore,  we focus on the LRT.

In the following, we denote $\bbX^{(1)}_{N_1}=(x_{ij}^{(1)})_{p\times N_1}$, where $\{x_{ij}^{(1)}\}$ are independent and identically distributed (i.i.d.) random variables with mean zero and variance one. Similarly $\bbX^{(2)}_{N_2}=(x_{ij}^{(2)})_{p\times N_2}$,  which is independent with $\bbX^{(1)}_{N_1}$,  constitutes another  i.i.d. sample with mean zero and variance one. For $l=1,2$,  we assume that  the $N_l$ observations $\bbz_j^{(l)}$ satisfy the {\it linear transformation model }
\begin{align}\label{ltm}
\bbz_j^{(l)}=\Sigma_l^{1/2}\bbx^{(l)}_{j}+\bm\mu_l,
\end{align}
where $\bbx_j^{(l)}$ is the $j$th column of $\bbX^{(l)}_{N_l}$ and $\bm\mu_l $ and $\Sigma_l$  are the population mean vectors and  covariance matrices   of $\{\bbz_j^{(l)}\}$, respectively.   Here,  $\Sigma_l^{1/2}$  can be chosen as any square root of matrix  $\Sigma_l$.  The linear transformation model covers the case where  the samples  $\bbz_i^{(l)}$ are  normally distributed, although  we are not  restricted here.  Denote $n_l=N_l-1$, $\bar\bbz^{(l)}=\frac{1}{N_l}\sum_{i=1}^{N_l}\bbz_i^{(l)}$,
$\mathring{\bbz}_i^{(l)}={\bbz}_i-\bar\bbz^{(l)}$ and $\bbS^{\bbz}_l:=\bbS_{n_l}^{(l)}=\frac{1}{n_l}\sum_{i=1}^{{N_l}}\mathring{\bbz}_i^{(l)}(\mathring{\bbz}_i^{(l)})'$.
We recall the famous Bartlett corrected LRT  statistic $L$ proposed  by \cite{Bartlett37P}, which is given by,
\begin{equation*}
L=\frac{2}{n_1+n_2}\log\(\frac{|\bbS^{\bbz}_1|^{\frac {n_1}{2}}\cdot|\bbS^{\bbz}_2|^{\frac{n_2}{2}}}{|{c_1}\bbS^{\bbz}_1+{c_2}\bbS^{\bbz}_2|^{\frac{n_1+n_2}{2}}}\).
\end{equation*}
%Because the LRT derived under the assumption of normality usually works well when the samples are drawn from non-normal populations, the same expressions are also called the LRT, although they are not.
Here, and in the following, we denote $c_1=n_1/(n_1+n_2)$ and $c_2=n_2/(n_1+n_2)$.
Moreover,
under the null hypothesis of  (\ref{hypo})
%{\color{blue}(all conclusions in this paper are basically derived under the null hypothesis)}  
and linear transformation model (\ref{ltm}), it  is not  difficult to rewrite  this
statistic as
\begin{align*}
L&
%\frac{2}{n_1+n_2}\log \(\frac{|\bbS^{\bbz}_1|^{\frac {N_1}{2}}\cdot|\frac{n_2}{N_2}\bbS^{\bbz}_2|^{\frac{N_2}{2}}}{|\frac{n_1}{N_1+N_2}\bbS^{\bbz}_1+\frac{n_2}
%	{N_1+N_2}\bbS^{\bbz}_2|^{\frac{N_1+N_2}{2}}}\)\\
=\frac{2}{n_1+n_2}\log \(\frac{|\bbS^{\bbx}_1|^{\frac {n_1}{2}}\cdot|\bbS^{\bbx}_2|^{\frac{n_2}{2}}}{|c_1\bbS^{\bbx}_1+c_2\bbS^{\bbx}_2|^{\frac{n_1+n_2}{2}}}\),
%=\frac{N_1}{2}\log(n_1\bbS_1(n_1\bbS_1+n_2\bbS_2)^{-1})+\frac{N_2}{2}\log(n_2\bbS_2(n_1\bbS_1+n_2\bbS_2)^{-1})
\end{align*}
%and
%\begin{align*}
%\tilde L=\log|n_1\bbS^{\bbx}_1(n_1\bbS^{\bbx}_1+n_2\bbS^{\bbx}_2)^{-1}|,
%\end{align*}
where  $\bbS^{\bbx}_l=\frac{1}{n_l}\sum_{i=1}^{{N_l}}(\bbx_i^{(l)}-\frac{1}{N_l}\sum_{i=1}^{N_l}
\bbx_i^{(l)})(\bbx_i^{(l)}-\frac{1}{N_1}\sum_{i=1}^{N_l}\bbx_i^{(l)})'$, $l=1,2$. Thus we know  that   $L$ is independent with  the population means $\bm\mu_l$ and covariance matrices $\Sigma_l$ under the null hypothesis $H_0$. In the following, we  drop the superscripts of $\bbS$ and
denote   $\bbS_1:=\bbS^{\bbx}_1$ and  $\bbS_2:=\bbS^{\bbx}_2$ for simplicity.
In addition, by simple calculation, we can rewrite
\begin{align*}
L
% =&\frac{N_1}{N_1+N_2}\log(n_1\bbS_1(n_1\bbS_1+n_2\bbS_2)^{-1})+\frac{N_2}{N_1+N_2}\log(n_2\bbS_2(n_1\bbS_1+n_2\bbS_2)^{-1})\\
= c_1\log(c_1^{-1}|\bbB_n|)+c_2\log(c_2^{-1}|\bbI_p-\bbB_n|),
\end{align*}
where
$$\bbB_n=n_1\bbS_1(n_1\bbS_1+n_2\bbS_2)^{-1}.$$

We now analyze $L$.
If $p>n_1$ or $p>n_2$, then $L$  is undefined because: (1) if $p\geq n_1+n_2$, then matrix $n_1\bbS_1+n_2\bbS_2$ is singular, which makes  the inverse of matrix  $n_1\bbS_1+n_2\bbS_2$ undefined; and (2) if $p< n_1+n_2$, i.e., the inverse of  $n_1\bbS_1+n_2\bbS_2$ is well-defined almost surely (with the fourth moments finite assumption),  but  $p>n_1$ or $p>n_2$,   then at least one of the determinants of
$\bbB_n$ or $\bbI_p-\bbB_n$  is zero, which makes the logarithm functions undefined.  However, from random matrix theory (RMT), we know that   if the fourth moment of $x_{ij}^{(l)}$ exists,  matrix $\bbB_n$ almost certainly has $p-n_1$ zero eigenvalues and $p-n_2$ one eigenvalues according to the condition  $p>n_1$ and $p>n_2$, respectively (see \cite{BaiH15C}).
Therefore,  we can naturally redefine the LRT $L$  by  restricting the  non-zero and non-one eigenvalues of $\bbB_n$, i.e.,
\begin{align}\label{lrts}
\mathscr L=\sum_{ \la_i^{\bbB_n}\in(0,1)}[c_1\log \la_i^{\bbB_n}+c_2\log (1-\la_i^{\bbB_n})],
%~~\mbox{and}~~
%\tilde L=\sum_{ \la_i^{\bbB_n}\in(0,1)}\log \la_i^{\bbB_n}
\end{align}
where $\lambda_i^{\bbB_n}$  denotes the i-th smallest eigenvalue of $\bbB_n$.
Therefore, we only need to obtain the asymptotic distributions of the redefined LRT in \eqref{lrts}, which is  addressed in the next section.

This paper also considers the  test statistic $\tilde L$,
\begin{equation*}
\tilde L=\log|n_1\bbS^{\bbz}_1(n_1\bbS^{\bbz}_1+n_2\bbS^{\bbz}_2)^{-1}|,
\end{equation*}
which was proposed by  \cite{PillaiJ67P,PillaiJ68P} and can be viewed as a lite LRT. Similar to $\mL$, we redefine $\tilde L$ by
\begin{align}\label{lrtts}
\tilde\mL=\sum_{ \la_i^{\bbB_n}\in(0,1)}\log \la_i^{\bbB_n}.
%~~\mbox{and}~~
%\tilde L=\sum_{ \la_i^{\bbB_n}\in(0,1)}\log \la_i^{\bbB_n}
\end{align}
It is obvious that $\tilde\mL$ is  monotone for matrix $\bbB_n$; thus, it should be more powerful than $\mL$.

The rest of this paper is organized as follows. The main results  are presented  in Section 2,  including the asymptotic normality of $\mL$ and $\tilde\mL$ and their optimal properties, which  stand for  the modification of LRT cannot be improved anymore in our model settings. In Section 3, we present the simulation results for the proposed statistics  by comparison with that proposed by \cite{LiC12T} and \cite{CaiL13T}. In Section 5, we introduce a real data application to demonstrate the application of the proposed tests.
All technical details are relegated to the appendix.
We note that for the high-dimensional testing problem \eqref{hypo}, the exact distribution of the test statistic  is difficult to obtain when the  distributions  are  free.

\section{Main results}
In this section, we give the asymptotic distributions of  the redefined LRT in \eqref{lrts} and \eqref{lrtts}.  For the application, we also present consistent estimators for the fourth moments of the samples under the null hypothesis.
Before presenting  the main results, we give some notation and the optimal assumptions. In the following,  we denote   the indicator function by $\delta_{(\cdot)}$, the natural logarithm function by $\log(\cdot)$,  convergence  in distribution by $\stackrel{D}{\to}$,  and
\begin{gather*}
y_1:= y_{n1}= p/n_1,\quad y_2:= y_{n2}= p/n_2,\quad h:= h_n= \sqrt{y_1+y_2-y_1y_2},\\
c_1:= n_1/(n_1+n_2)=y_2/(y_1+y_2), \quad c_1:= n_2/(n_1+n_2)=y_1/(y_1+y_2),\\
l(y_1,y_2)=\log({h^{\frac{2c_1h^2}{y_{1}y_{2}}}})\delta_{y_1>1}-\log({y_1^\frac{c_1(1+y_2)}{y_2}y_2^{\frac{c_1(1-y_1)}{y_1}}}
)\delta_{y_1>1},\\
u(y_1,y_2)=\log(\frac{y_1^{c_1}}{h^{c_1}})\delta_{y_1>1},\quad v(y_1,y_2)=\log\left(\frac{y_1^{2c_1}}{h^{2c_1(c_1+2c_2)}}\right)\delta_{y_1>1},\\
\Psi(y_1,y_2)=c_2y_1^2[y_2^4\delta_{y_2<1}+h^2(2y_2^2-h^2)\delta_{y_2>1}]\qquad\qquad\qquad\\
\qquad\qquad\qquad-c_1y_2^2[y_1^3(y_1+2y_2)\delta_{y_1<1}+h^2(y_1+y_2+y_1y_2)\delta_{y_1>1}].
\end{gather*}
We set two assumptions of the sample that will be shown to be optimal for the proposed test statistics.
\begin{itemize}
	\item (Moments Assumption:) $\E x_{11}^{(1)}=\E x_{11}^{(2)}=0, \E(x_{11}^{(1)})^2=\E(x_{11}^{(2)})^2=1$, $\E(x_{11}^{(1)})^4=\Delta_1+3<\infty$ and  $\E(x_{11}^{(2)})^4=\Delta_2+3<\infty$;
	\item (Dimensions Assumption:) $y_1\neq1$, $ y_2\neq1$ and
	$ p/
	{(n_1+n_2)}<1$.
\end{itemize}
We are now ready to present the main results of this paper.
\begin{theorem} \label{th1}
	In addition to  the Moments Assumption and Dimensions Assumption, we assume that as $\min\{p,n_1,n_2\}\to\infty$,  $\lim y_1\not\in\{0,1\}$, $ \lim y_2\not\in\{0,1\}$ and
	$ \lim p/
	{(n_1+n_2)}\in(0,1)$.  Then under the null hypothesis, we have
	\begin{align}\label{TT}
	T:=	\frac{\mL-p\ell_n-\mu_n}{\nu_n}\stackrel{D}{\to}N(0,1),%~~~~\mbox{~~and~~}%\tilde T:=\frac{\tilde L-p\tilde l_n-\tilde \mu}{\tilde v}\stackrel{D}{\to}N(0,1),
	\end{align}
	where
	\begin{align*}
	\ell_n=&\log\left(\frac{y_1^{c_2} y_2^{c_1}h^{\frac{2h^2}{y_1y_2}}}{(y_1+y_2)^{\frac{(y_1+y_2)}{y_1y_2}}|1-y_1|^{\frac{c_1|1-y_1|}{y_1}}|1-y_2|^{\frac{c_2|1-y_2|}{y_2}}}
	\right)-l(y_1,y_2)-l(y_2,y_1),
	%&-\log\left[\frac{h^{\frac{2c_1h^2}{y_{1}y_{2}}}}{y_1^\frac{c_1(1+y_2)}{y_2}y_2^{\frac{c_1(1-y_1)}{y_1}}}
	%\right]\delta_{y_1>1}-\log\left[\frac{h^{\frac{2c_2h^2}{y_{1}y_{2}}}}{y_1^{\frac{c_2(1-y_2)}{y_2}}y_2^\frac{c_2(1+y_1)}{y_1}}
	%\right]\delta_{y_2>1}
	\end{align*}
	%
	%
	%\[l_n=\left\{\begin{array}{ll}
	%p\log[(y_{1}+y_{2})^{-\frac{(y_{1}+y_{2})}{y_{1}y_{2}}}
	%\times(y_{1}-1)^{\frac{c_1(1-y_{1})}{y_{1}}}\times(y_{2}-1)^{\frac{c_2(1-y_{2})}{y_{2}}}
	%\times y_{1}^{\frac{(c_1y_{2}+1)}{y_{2}}}\times y_{2}^{\frac{(c_2y_{1}+1)}{y_{1}}}
	%]\\
	%\ \hspace{340pt}\text{,$\gamma_1>1$, $\gamma_2>1$}\\
	%p\log[(y_{1}+y_{2})^{-\frac{(y_{1}+y_{2})}{y_{1}y_{2}}}\times(y_{1}-1)^{\frac{c_1(1-y_{1})}{y_{1}}}
	%\times(1-y_{2})^{\frac{c_2(y_{2}-1)}{y_{2}}}
	%\times y_{1}^{\frac{(y_{2}+c_1)}{y_{2}}}
	% \times y_{2}^{\frac{c_1}{y_{1}}}\times (h^2)^{\frac{c_2h^2}{y_{1}y_{2}}}]\\
	% \ \hspace{340pt}\text{,$\gamma_1>1$, $\gamma_2<1$}\\
	% p\log[(y_{1}+y_{2})^{-\frac{(y_{1}+y_{2})}{y_{1}y_{2}}}
	% \times (1-y_{1})^{\frac{c_1(y_{1}-1)}{y_{1}}}
	% \times(y_{2}-1)^{\frac{c_2(1-y_{2})}{y_{2}}}
	% \times y_{1}^{\frac{c_2}{y_{2}}}
	% \times y_{2}^{\frac{(c_2+y_{1})}{y_{1}}}
	% \times (h^2)^{\frac{c_1h^2}{y_{1}y_{2}}}
	% ]\\
	% \ \hspace{340pt}\text{,$\gamma_1<1$, $\gamma_2>1$}\\
	%p\log[(y_{1}+y_{2})^{-\frac{(y_{1}+y_{2})}{y_{1}y_{2}}}
	%\times(1-y_{1})^{\frac{c_1(y_{1}-1)}{y_{1}}}
	%\times(1-y_{2})^{\frac{c_2(y_{2}-1)}{y_{2}}}
	%\times y_{1}^{c_2}
	%\times y_{2}^{c_1}
	%\times(h^2)^{\frac{h^2}{y_{1}y_{2}}}]\\
	% \ \hspace{340pt}\text{,$\gamma_1<1$, $\gamma_2<1$}\\
	%\end{array}\right.\]
	%
	%
	%
	\begin{align*}
	\mu_n=&\log\left[\frac{{(y_1+y_2)^\frac{1}{2}|1-y_1|^\frac{c_1}{2}|1-y_2|^\frac{c_2}{2}}}{ h}\right]-u(y_1,y_2)-u(y_2,y_1)\\
	&+\frac{\Delta_1\Psi(y_1,y_2)}{2y_1y_2^2(y_1+y_2)^2}
	+\frac{\Delta_2\Psi(y_2,y_1)}{2y_2y_1^2(y_1+y_2)^2},
	%&+\frac{\Delta_2}{2y_2y_1^2(y_1+y_2)^2}\left[-c_2y_2^3y_1^2(y_2+2y_1)\delta_{y_2<1}
	%+c_1y_2^2y_1^4\delta_{y_1<1}-{c_2y_1^2h^2(y_1+y_2+y_1y_2)}
	%\delta_{y_2>1}+{c_1y_2^2h^2(2y_1^2-h^2)}\delta_{y_1>1}\right]
	\end{align*} and
	\begin{align*}
	\nu^2_n=&\log\frac{h^4}{|1-y_1|^{2c_1^2}|1-y_2|^{2c_2^2}(y_1+y_2)^2}+2[v(y_1,y_2)+v(y_2,y_1)+\log(h^{4c_1c_2})\delta_{y_1>1}\delta_{y_2>1}]\\
	&+\frac{(y_1\Delta_1+y_2\Delta_2)}
	{y_1^2y_2^2(y_1+y_2)^2}[(y_1-1)y_2^2\delta_{y_1>1}-(y_2-1)y_1^2\delta_{y_2>1}]^2.
	\end{align*}
\end{theorem}

%{\color{blue}The asymptotic distribution of the proposed test statistic $T$ under the null hypothesis is displayed through the Theorem \ref{th1}. We can get an accordingly critical region $W=\{|T|>C\}$ based on a given significance level $0.05$, where $C$ satisfies $\mathbb{P}(|T|>C|H_0)=0.05$. One should reject the null hypothesis when observed value is in the critical region. }
%
%{\color{blue}Theorem \ref{th1} follows from first formulation of (1.2) in \cite{BaiH15C} and the proof of the Theorem \ref{th1} is supplied as supplementary material for interested readers, see \cite{ZhangH17M}.}

\begin{figure}[htbp]
	\centering
	\subfigure[3D figure of the mean $\mu_n$.]{
		\label{fig:subfig:a} %% label for first subfigure
		\includegraphics[scale=0.45]{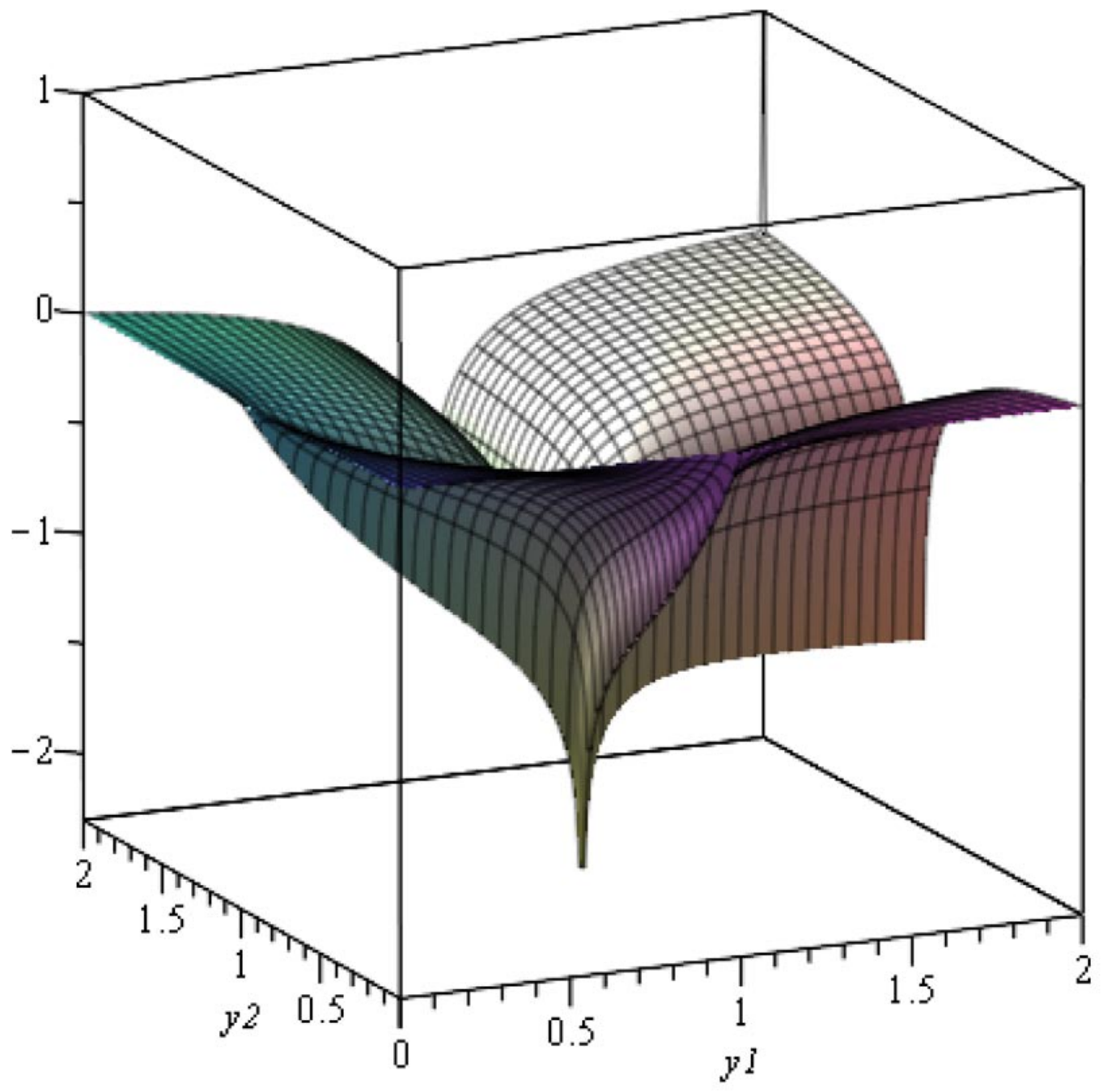}}
	%	\hspace{1in}
	\subfigure[3D figure of the variance $\nu_n^2$.]{
		\label{fig:subfig:b} %% label for second subfigure
		\includegraphics[scale=0.41]{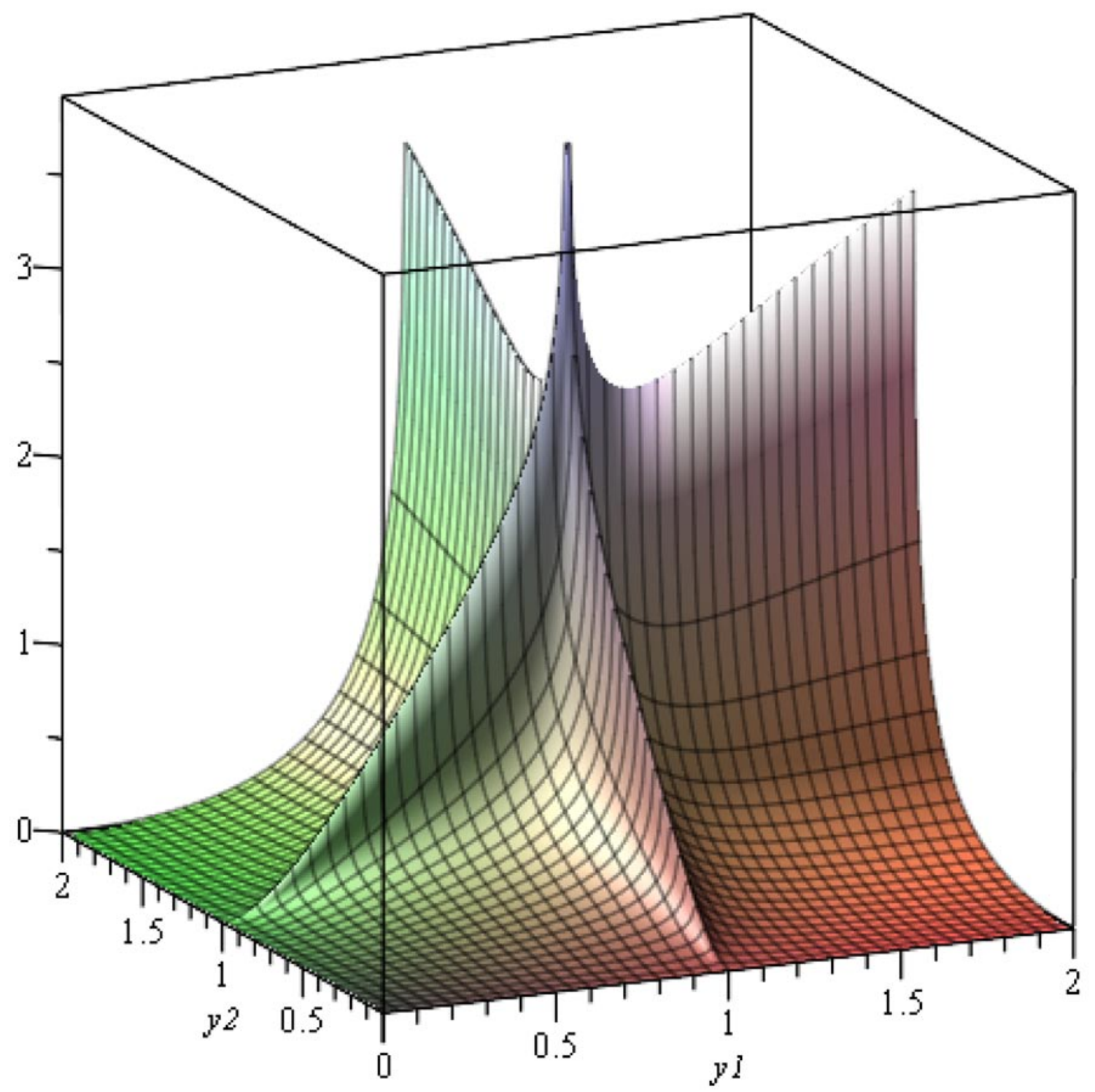}}
	\caption{ This figure was made using Maple software with $y_1\in(0,2)$ and $y_2\in(0,2)$. The vertical axes present the values of $\mu_n$ and  $\nu_n$, respectively. }\label{figclt}
\end{figure}

According to the above theorem, we can easily conclude the following corollary under normal circumstances:

\begin{corollary} \label{cor1}If $z_j^{(l)}$, $j=1\dots,n_l$, $l=1,2$ are normally distributed, then \eqref{TT} in Theorem \ref{th1}  reduces to:
	\begin{align*}
	T:=	\frac{\mL-p\ell_n-\mu_n}{\nu_n}\stackrel{D}{\to}N(0,1),%~~~~\mbox{~~and~~}%\tilde T:=\frac{\tilde L-p\tilde l_n-\tilde \mu}{\tilde v}\stackrel{D}{\to}N(0,1),
	\end{align*}
	where
	\begin{align*}
	\ell_n=&\log\left(\frac{y_1^{c_2} y_2^{c_1}h^{\frac{2h^2}{y_1y_2}}}{(y_1+y_2)^{\frac{(y_1+y_2)}{y_1y_2}}|1-y_1|^{\frac{c_1|1-y_1|}{y_1}}|1-y_2|^{\frac{c_2|1-y_2|}{y_2}}}
	\right)-l(y_1,y_2)-l(y_2,y_1),
	\end{align*}
	\begin{align*}
	\mu_n=&\log\left[\frac{{(y_1+y_2)^\frac{1}{2}|1-y_1|^\frac{c_1}{2}|1-y_2|^\frac{c_2}{2}}}{ h}\right]-u(y_1,y_2)-u(y_2,y_1),
	\end{align*} and
	\begin{align*}
	\nu^2_n=&\log\frac{h^4}{|1-y_1|^{2c_1^2}|1-y_2|^{2c_2^2}(y_1+y_2)^2}+2[v(y_1,y_2)+v(y_2,y_1)+\log(h^{4c_1c_2})\delta_{y_1>1}\delta_{y_2>1}].
	\end{align*}
	
\end{corollary}
\begin{remark}
	If the Dimensions Assumption in Theorem \ref{th1} is satisfied, the limit condition that  $\lim y_1\not\in\{0,1\}$, $ \lim y_2\not\in\{0,1\}$ and $\lim p/(n_1+n_2)\in(0,1)$ could be considered to  hold,  because there is no information for the convergence of its dimensions and   size for any dataset. In addition, \cite{JiangY13C}  showed that if $p/n_l\to1$, Theorem \ref{th1} also holds for normally distributed data.
	%The limit condition here guarantee s $\mu_n$ and $\nu_n$ are not to too small and too large when the dimensions are large.
	However, if $y_l$ is near 1, then the variance $\nu_n$ will be   large  and the  LRT will become unstable, see Figure \ref{figclt} for illustration.
\end{remark}
\begin{remark}
	When $y_1<1$ and $y_2<1$, Corollary \ref{cor1}  recovers Theorem 4.1 in \citep{BaiJ09C}  directly.
	
\end{remark}
The optimality of the Dimensions Assumption is clear because of the definitions of ${\ell_n}$, $\mu_n$ and $\nu_n$. For the Moments Assumption, we only need to consider the fourth moments of the sample. From the variance $\nu_n$ in Theorem \ref{th1}, we know that its fourth-moment  term cannot be removed except $y_1<1$ and $y_2<1$. However, if $y_1<1$ and $y_2<1$, it is not difficult to obtain that \begin{align*}
\Psi(y_2,y_1)=\Psi(y_1,y_2)=-y_1^3y_2^3,
\end{align*}
which implies that the fourth-moment term of $\mu_n$ cannot be removed. Therefore, we conclude that the existence of the fourth moment of the sample is necessary  for the modified LRT statistic $\mL$. However, when $y_1$ or $y_2$ is close to 1, the variance $\nu_n$ will increase rapidly, resulting in poor power.   For illustration, we present two 3D figures of $\mu_n$ and $\nu_n^2$ with $\Delta_1=\Delta_2=0$ in Figure \ref{figclt}.

Now, we give the asymptomatic distribution of $\tilde\mL$.
\begin{theorem} \label{th2}
	In addition to the Moments Assumption and Dimensions Assumption, we assume that as $\min\{p,n_1,n_2\}\to\infty$,  $\lim y_1\neq1$, $ \lim y_2\neq1$ and
	$ \lim p/
	{(n_1+n_2)}<1$.  Then under the null hypothesis, we have
	\begin{align*}
	\tilde T:=	\frac{\tilde \mL-p\tilde \ell_n-\tilde \mu_n}{\tilde \nu_n}\stackrel{D}{\to}N(0,1),
	\end{align*}
	where
	\begin{align*}
	\tilde \ell_n=&\log\left(\frac{ y_2h^{\frac{2h^2}{y_1y_2}}}{(y_1+y_2)^{\frac{(y_1+y_2)}{y_1y_2}}|1-y_1|^{\frac{|1-y_1|}{y_1}}}
	\right)
	-\log\(\frac{{h^{\frac{2h^2}{y_1y_2}}}}{{y_1^\frac{1+y_2}{y_2}y_2^{\frac{1-y_1}{y_1}}}}\)\delta_{y_1>1}
	\end{align*}
	\begin{align*}
	&\tilde \mu_n=\log\left[\frac{{(y_1+y_2)^\frac{1}{2}|1-y_1|^\frac{1}{2}}}{ h}\right]-\log(\frac{y_1}{h})\delta_{y_1>1}\\
	&	-\frac{\Delta_1[y_1^3(y_1+2y_2)\delta_{y_1<1}+h^2(y_1+y_2+y_1y_2)\delta_{y_1>1}]}{2y_1(y_1+y_2)^2}
	+\frac{\Delta_2[y_1^4y_2\delta_{y_1<1}+h^2y_2(2y_1^2-h^2)\delta_{y_1>1}]}{2y_1^2(y_1+y_2)^2},
	\end{align*}
	and
	\begin{align*}
	\tilde \nu^2_n=&2\log\frac{h^2}{|1-y_1|(y_1+y_2)}+2\log\(\frac{y_1^2}{h^2}\)\delta_{y_1>1}+\frac{(y_1\Delta_1+y_2\Delta_2)}
	{y_1^2(y_1+y_2)^2}[y_1^4\delta_{y_1<1}+h^4\delta_{y_1>1}].
	\end{align*}
\end{theorem}
\begin{remark}
	Notice that the asymptotic distributions  in Theorem \ref{th1} and Theorem \ref{th2} are obtained under the null hypothesis, which can only guarantee the Type I errors. For the powers, that is  under the alternative hypothesis  $\Sigma_1\neq \Sigma_2$, the asymptotic distributions of statistics $\mL$ and $\tilde\mL$ will depend on the eigenvalues of $\Sigma_1\Sigma_2^{-1}$. In this case, if the dimension $p$ is smaller than either of the two  sample sizes, the power functions  for  $T$ and $\tilde{T}$ can be obtained by the CLT  of  the 
	general Fisher matrices which is  derived by \cite{ZhengB17C}. However, if $p$ is bigger than  both of the  two  sample sizes, because of the lack of  theoretical   results about the general Beta matrix $n_1\bbS_1(n_1\bbS_1+n_2\bbA^{1/2}\bbS_2\bbA^{1/2})^{-1}$,  the asymptotic distributions of statistics $\mL$ and $\tilde\mL$ are also open problems and will be left for our future work.  Here $\bbA$ is any non-random symmetric matrix.
\end{remark}

If $\Delta_l\neq0$, or more specifically, if  the samples are not normally distributed,  the estimates for  $\Delta_l$ are necessary for the test application. Thus, we obtain   their  consistent estimators
using the method of moments and random  matrix theory.
{
Let \begin{align}
\hat\Delta_1=(1-y)^2\frac{\sum_{j=1}^{N_1}[(\bbz^{(1)}_j-\bar{\bbz}^{(1)})'(c_{11}\bbS^
{\bbz}_{1j}+c_{12}\bbS_2^{\bbz})^{-1}(\bbz^{(1)}_j-\bar{\bbz}^{(1)})-\frac{p}{1-y}]^2}{pN_1}
-\frac{2}{1-y}\label{delta1}\\
\hat\Delta_2=(1-y)^2\frac{\sum_{j=1}^{N_2}[(\bbz^{(2)}_j-\bar{\bbz}^{(2)})'(c_{21}\bbS^
{\bbz}_{1}+c_{22}\bbS_{2j}^{\bbz})^{-1}(\bbz^{(2)}_j-\bar{\bbz}^{(2)})-\frac{p}{1-y}]^2}{pN_2}
-\frac{2}{1-y},\label{delta2}
\end{align}
where $y=\frac{p}{n_1+n_2-1}$, $c_{11}=\frac{n_1-1}{n_1+n_2-1}$, $c_{12}=\frac{n_2}{n_1+n_2-1}$,  $c_{21}=\frac{n_1}{n_1+n_2-1}$, $c_{22}=\frac{n_2-1}{n_1+n_2-1}$ and $\bbS^{\bbz}_{lj}$ is the sample covariance matrix by removing the vector $\bbz_{j}^{(l)}$ from the $l$-th sample, $l=1,2$.}
{
\begin{theorem}\label{th3}
	Under the same assumptions of Theorem \ref{th1} and under the null hypothesis, we have the estimators $\hat\Delta_l$, $l=1,2$,   defined in \eqref{delta1}  and \eqref{delta2} are weakly consistent and asymptotically unbiased.
\end{theorem}
\begin{remark}
Actually we use  	$(c_{11}\bbS^{\bbz}_{1j}+c_{12}\bbS_2^{\bbz})^{-1}$ here is to avoid the case where the sample covariance matrices when the sample sizes are not larger than  the dimension. Otherwise,
if the dimension $p$  is smaller than some sample size, say, $p<n_1$, we can estimate $\Delta_1$ by replacing the term $(c_{11}\bbS^{\bbz}_{1j}+c_{12}\bbS_2^{\bbz})^{-1}$ and $y=\frac{p}{n_1+n_2-1}$ in \eqref{delta1} by  $(\bbS^{\bbz}_{1j})^{-1}$ and  $y=\frac{p}{n_1-1} $ respectively.
\end{remark}}
	The proofs of Theorem \ref{th1}, Theorem \ref{th2} and Theorem \ref{th3} are  given in the appendix.

\section{Results of the simulation}

In this section, we compare the performance of the statistics proposed in \cite{LiC12T} and \cite{CaiL13T} and our modified LRTs $\mL$ and $\tilde{\mL}$  under various settings of sample size and dimensionality.  The classical LRT  statistic in \citep{Wilks46S} was shown to have poor performance for (\ref{hypo})
by \cite{BaiJ09C};  thus, it will not be considered in this section. Without loss of generality, we assume $\bm\mu_l=0$ and set
$$\Sigma_1=(1+a/n_1)\Sigma_2,$$
where  $a$ is a constant.
The samples are drawn from the following  distributions:
\begin{itemize}
	\item[]Case 1: $\bbx^{(1)}$ and $\bbx^{(2)}$ are both  standard normal distributed and $\Sigma_2=\bbI_p$;
	\item[]  Case 2: $\bbx^{(1)}$ and $\bbx^{(2)}$  are  from the uniform distribution  $U(-\sqrt{3},\sqrt{3})$ and $\Sigma_2=\bbI_p$;
	\item[]  Case 3: $\bbx^{(1)}$ and $\bbx^{(2)}$  are  from the uniform distribution  $U(-\sqrt{3},\sqrt{3})$ and $$\Sigma_2=Diag(p^2,1,\dots, 1);$$
	\item[] Case 4: $\bbx^{(1)}$ and $\bbx^{(2)}$  are  from the uniform distribution  $U(-\sqrt{3},\sqrt{3})$ and $$\Sigma_2=(0.5\bbI_p+0.5\textbf{1}_p \textbf{1}_p ').$$
\end{itemize}
Here, $\textbf{1}_p$ represents a $p$-dimensional vector with all entries 1. The
results are obtained based on 10,000 replicates. In the tables,
$T$ and $\tilde{T}$
denote the proposed modified LRTs, $T_{lc}$ denotes the nonparametric test of  \cite{LiC12T} and $T_{clx}$ denotes the maximum element  test of \cite{CaiL13T}.

\begin{figure}[htbp]
	\begin{center}
		\includegraphics[width=14cm]{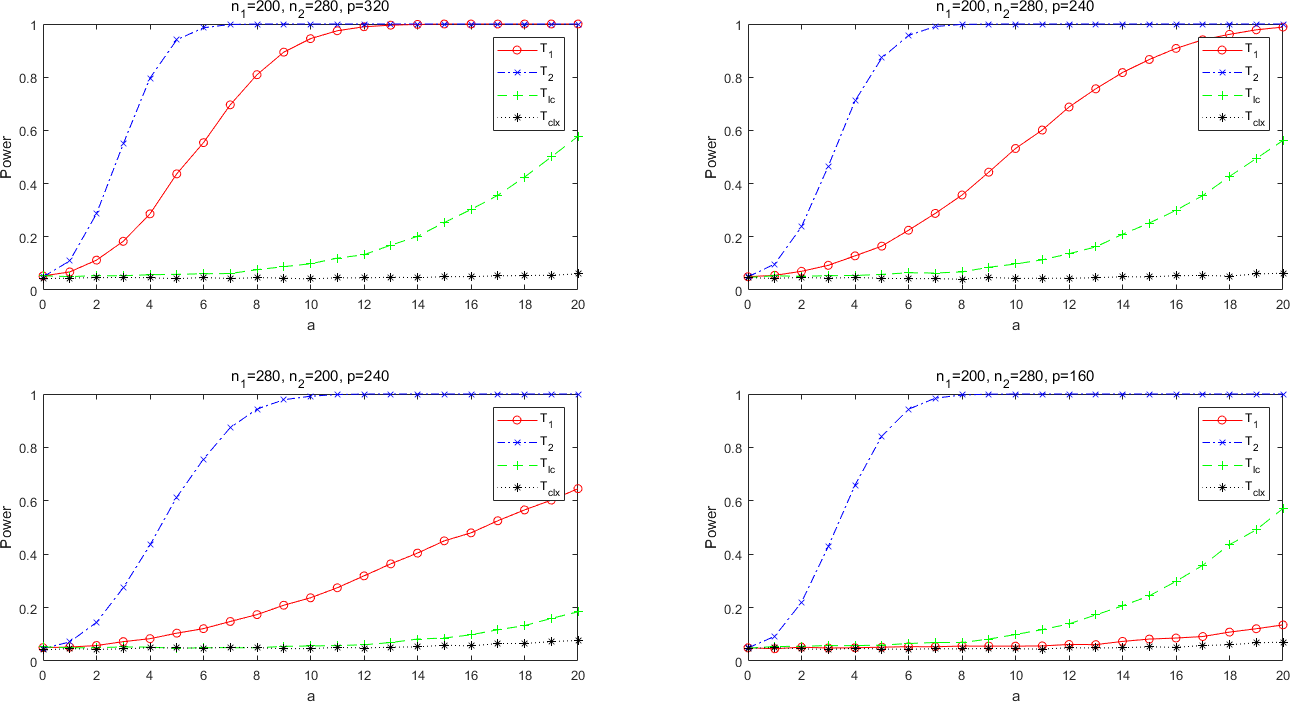}
	\end{center}
	\caption{Graphs of the  divergence of the four powers  in Case 1.}\label{fig2}
\end{figure}

\begin{figure}[htbp]
	\begin{center}
		\includegraphics[width=14cm]{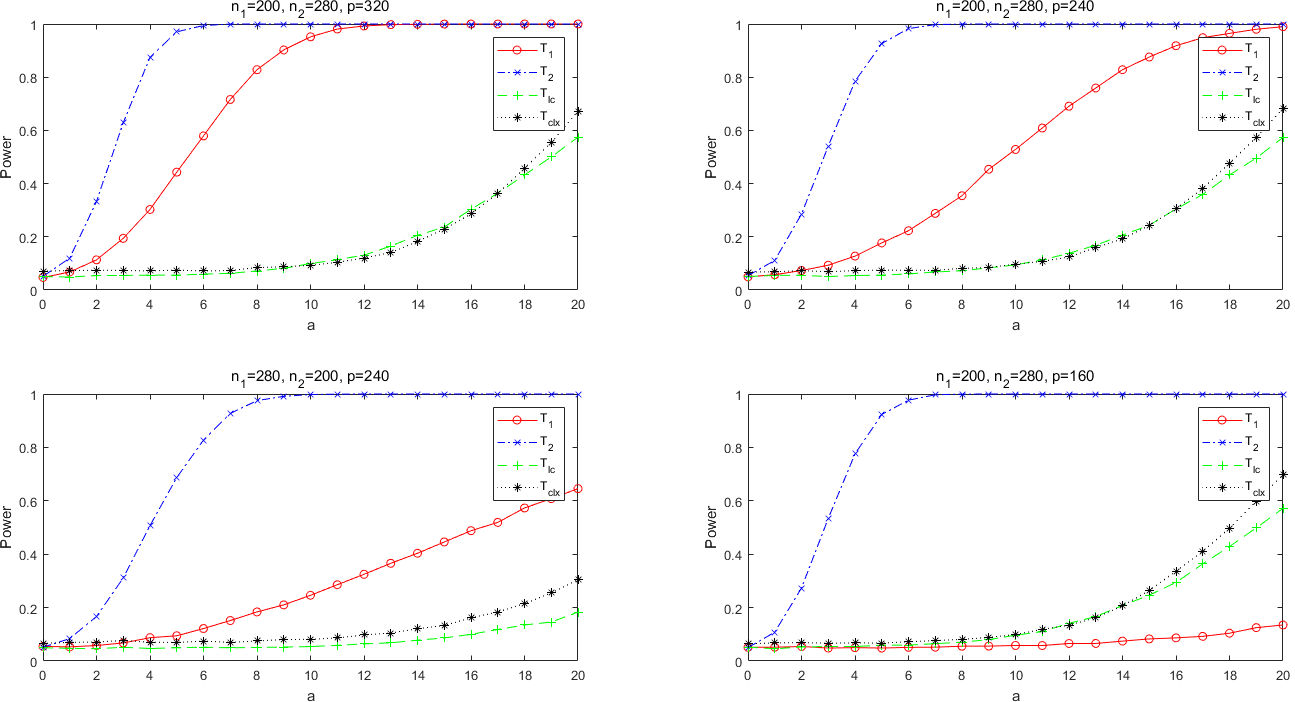}
	\end{center}
	\caption{Graphs of the  divergence of the four powers  in Case 2.}\label{fig3}
\end{figure}

\begin{figure}[htbp]
	\begin{center}
		\includegraphics[width=14cm]{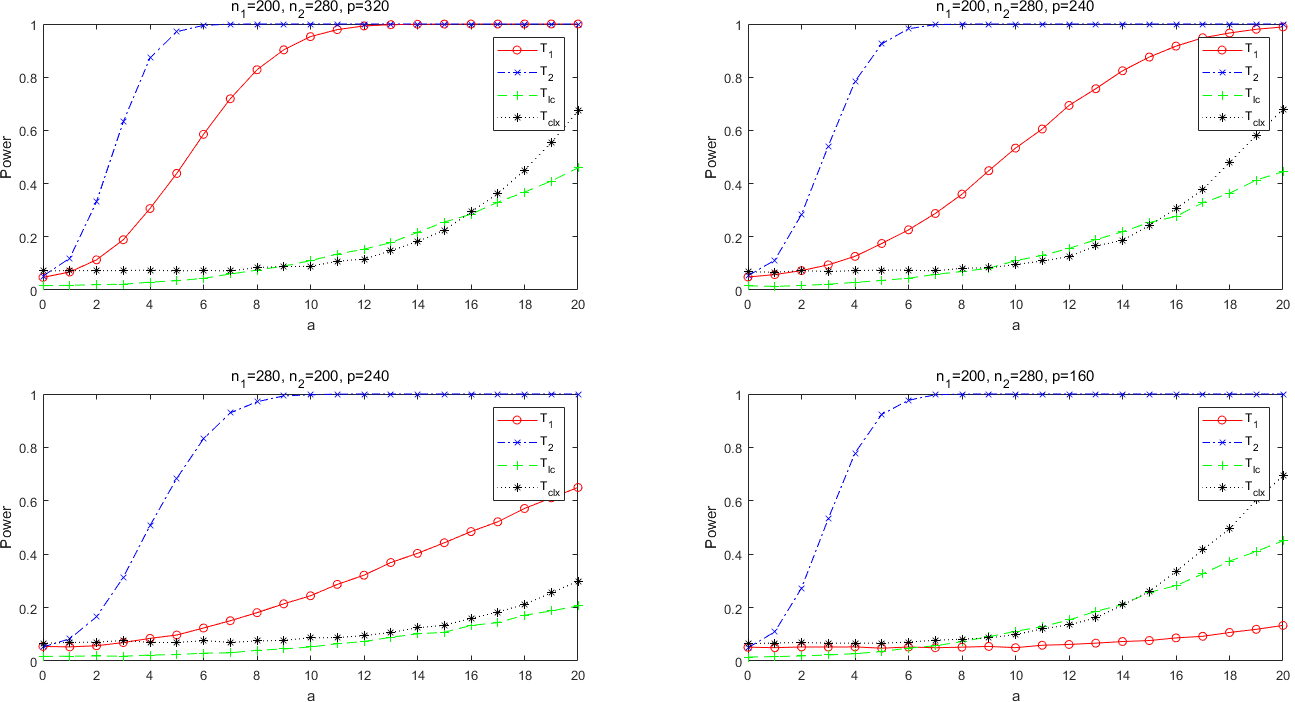}
	\end{center}
	\caption{Graphs of the  divergence of the four powers  in Case 3.}\label{fig4}
\end{figure}

\begin{figure}[htbp]
	\begin{center}
		\includegraphics[width=14cm]{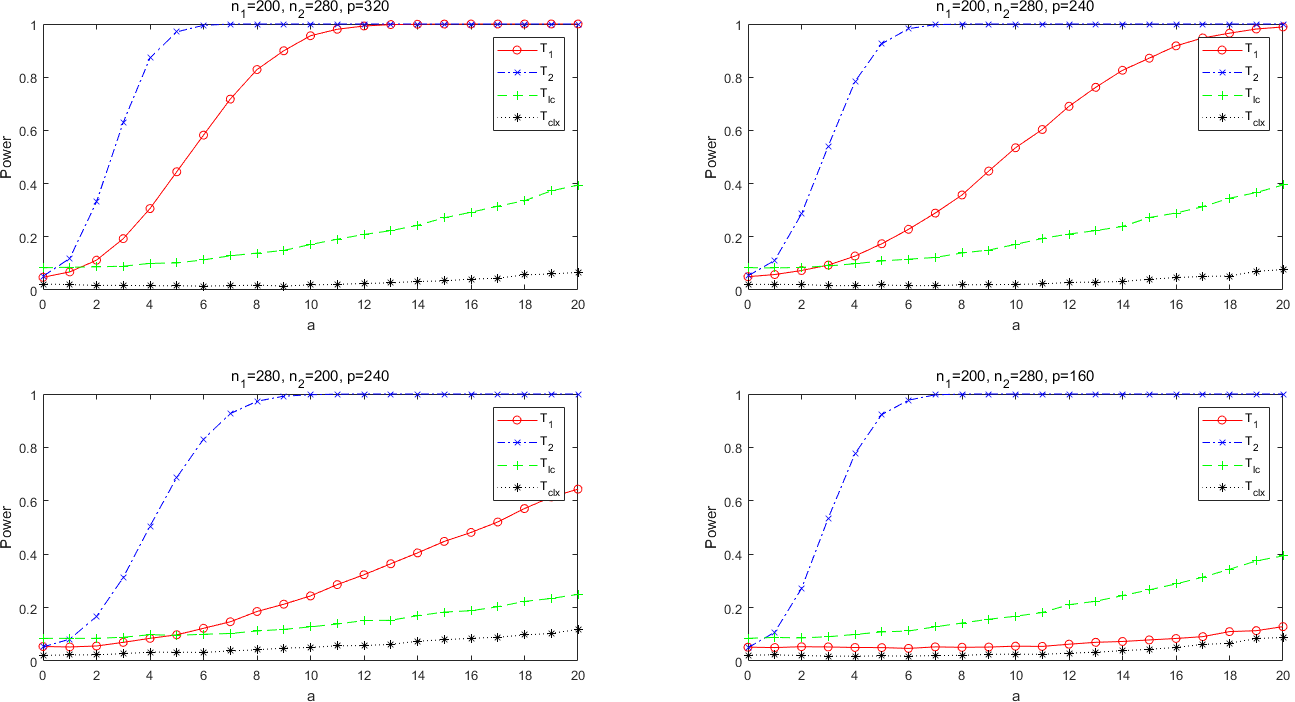}
	\end{center}
	\caption{Graphs of the  divergence of the four powers  in Case 4.}\label{fig5}
\end{figure}

In the first part of this section, we report the results by assuming the forth moments of the  $\bbx^{(1)}$ and $\bbx^{(2)}$ are known.
Tables \ref{tab1}- \ref{tab4} present the empirical sizes and empirical powers of Cases 1- 4. Additionally, we provide four figures (Figures \ref{fig2} - \ref{fig5}) to show the divergence of powers of the four test statistics as the parameter $a$ increases.
The results indicate that the sizes $T$ and $\tilde T$ perform quite well for all cases. However, in Case 3 and  Case 4, the sizes of  $T_{lc}$ and $T_{clx}$  are not  accurate, which reflects the fact that the null distributions of the test statistics $T_{lc}$ and $T_{clx}$
are not  well approximated by their asymptotic distributions in these cases.  For the empirical power,  we can  conclude  that $T$ and $\tilde T$  are more sensitive than $T_{lc}$  and $T_{clx}$ when at least one of the sample sizes is smaller than the dimensions. Otherwise, the LRT $T$ does not perform as well when the dimensions are large. However, the lite LRT $\tilde T$ is always powerful because of its monotonicity  for matrix $\bbB_n$, which  coincides with  our intuition.

{
Next we will show the performance of the estimator we proposed   in Theorem \ref{th3}. Here we have to explain  the reason that why we did not use the estimators in last simulation results. That is because
our estimator is based on moment method, and need  $n=n_1+n_2$ times loop and inverse process for one replication, and we need 10,000 replications for one result, that  makes the running times to be terrible.
In the simulation,
we set $x_{ij}$  be  standard normal distributed and uniform distributed on ($-\sqrt3,\sqrt3$) to estimate the $\Delta_1$ respectively. Under each circumstance we repeat  10,000 times and the results are reported at Tables \ref{Del1} and \ref{Del2}. From the numerical results, the performance of the estimator  is remarkable, especially when the sample size is large. Therefore, we believe that  the proposed modified LRTs must be also perform good at the null hypothesis when  using the estimators instead of their true values. But, under the alternative, we can not make any arbitrary decision right now because of the less  theoretical results about the general Beta  matrix. }

\begin{table}[htbp]
\scriptsize
	\begin{tabular}{ccccccccccccc}
		%		\multicolumn{13}{c}
		%		{Case1} \\
		\hline
		\multirow{3}{6.3em}{{\centering ($n_1,n_2,p$)\\$y_1>1$,$y_2>1$}}&\multicolumn{3}{c}{(25,35,40)} &\multicolumn{3}{c}{(50,70,80)} &\multicolumn{3}{c}{(100,140,160)} &\multicolumn{3}{c}{(200,280,320)}\\ \cline{2-13}
		&size&\multicolumn{2}{c}{power}
		&size&\multicolumn{2}{c}{power}&size&\multicolumn{2}{c}{power}&size&\multicolumn{2}{c}
		{power}\\\cline{2-13}
		&$a$=0&$a$=10&$a$=20&$a$=0 &$a$=10&$a$=20&$a$=0&$a$=10&$a$=20&$a$=0&$a$=10&$a$=20 \\ \hline
		$T  $                   &0.055 &0.951&1     & 0.053&0.951&  1   &0.054&0.951&   1  &0.051&0.946 & 1     \\ \hline
		$\tilde{T}$              &0.056 &0.999&1     &0.048 &1    &1     &0.048&1    &1     &0.049&1     & 1   \\ \hline
		$T_{lc}$                &0.074 &0.804&1     &0.057 &0.489&0.999 &0.054&0.208&0.954 &0.054&0.093 &0.571  \\ \hline
		$T_{clx}$                     &0.082 &0.152&0.591 &0.057 &0.071&0.343 &0.048&0.049&0.120 &0.042&0.046 &0.062  \\ \hline
		
		\multirow{3}{6.3em}{({\centering$n_1,n_2,p$)\\$y_1>1$,$y_2<1$}}&\multicolumn{3}{c}{(25,35,30)} &\multicolumn{3}{c}{(50,70,60)} &\multicolumn{3}{c}{(100,140,120)} &\multicolumn{3}{c}{(200,280,240)}\\ \cline{2-13}
		&size&\multicolumn{2}{c}{power}
		&size&\multicolumn{2}{c}{power}&size&\multicolumn{2}{c}{power}&size&\multicolumn{2}{c}
		{power}\\\cline{2-13}
		&$a$=0&$a$=10&$a$=20&$a$=0 &$a$=10&$a$=20&$a$=0&$a$=10&$a$=20&$a$=0&$a$=10&$a$=20 \\ \hline
		$T  $                   &0.059&0.685&0.999 & 0.057&0.618&0.998 &0.054&0.558& 0.994&0.052&0.523 &0.987     \\ \hline
		$\tilde{T}$              &0.058&0.997&1     & 0.053&0.999&  1   &0.050&0.999&   1  &0.049& 1    & 1     \\ \hline
		$T_{lc}$                &0.063&0.793&1     & 0.055&0.494&0.999 &0.053&0.205&0.952 &0.050&0.098 &0.574  \\ \hline
		$T_{clx}$                     &0.081&0.154&0.621 &0.054 &0.076&0.380 &0.049&0.048&0.147 &0.045&0.046 &0.068  \\ \hline
		\multirow{3}{6.3em}{{\centering ($n_1,n_2,p$)\\$y_1<1$,$y_2>1$}}&\multicolumn{3}{c}{(35,25,30)} &\multicolumn{3}{c}{(70,50,60)} &\multicolumn{3}{c}{(140,100,120)} &\multicolumn{3}{c}{(280,200,240)}\\ \cline{2-13}
		&size&\multicolumn{2}{c}{power}
		&size&\multicolumn{2}{c}{power}&size&\multicolumn{2}{c}{power}&size&\multicolumn{2}{c}
		{power}\\\cline{2-13}
		&$a$=0  &$a$=10 &$a$=20 &$a$=0  &$a$=10 &$a$=20 &$a$=0  &$a$=10 &$a$=20 &$a$=0  &$a$=10 &$a$=20  \\ \hline
		$T$                   &0.060	&0.109	&0.075	&0.055	&0.167	&0.282	&0.050	&0.213	&0.515	&0.050	&0.248	&0.645    \\ \hline
		$\tilde{T}$              &0.056	&0.956	&1	    &0.052	&0.985	&1	    &0.052	&0.989	&1	    &0.046	&0.994	&1        \\ \hline
		$T_{lc}$                &0.068	&0.314	&0.992	&0.059	&0.143	&0.896	&0.054	&0.077	&0.489	&0.051	&0.054	&0.179    \\ \hline
		$T_{clx}$                     &0.077	&0.203	&0.677	&0.054	&0.102	&0.433	&0.049	&0.061	&0.189	&0.042	&0.051	&0.076   \\ \hline
		\multirow{3}{6.3em}{{\centering ($n_1,n_2,p$)\\$y_1<1$,$y_2<1$}}&\multicolumn{3}{c}{(25,35,20)} &\multicolumn{3}{c}{(50,70,40)} &\multicolumn{3}{c}{(100,140,80)} &\multicolumn{3}{c}{(200,280,160)}\\ \cline{2-13}
		&size&\multicolumn{2}{c}{power}
		&size&\multicolumn{2}{c}{power}&size&\multicolumn{2}{c}{power}&size&\multicolumn{2}{c}
		{power}\\\cline{2-13}
		&$a$=0&$a$=10&$a$=20&$a$=0 &$a$=10&$a$=20&$a$=0&$a$=10&$a$=20&$a$=0&$a$=10&$a$=20 \\ \hline
		$T _1$                   &0.060&0.251&0.960 & 0.056&0.119&0.715 &0.050&0.079&0.330 &0.049&0.053 & 0.139     \\ \hline
		$\tilde{T}$              &0.057&0.996&1     & 0.052&0.999&  1   &0.049&1    &   1  &0.050&0.999 & 1     \\ \hline
		$T_{lc}$                &0.069&0.783&1     & 0.059&0.476&0.999 &0.051&0.211&0.948 &0.049&0.095 &0.563  \\ \hline
		$T_{clx}$                     &0.076&0.169&0.679 &0.058 &0.083&0.438 &0.046&0.053&0.180 &0.047&0.050 &0.072  \\ \hline
		
	\end{tabular}
	\caption{Empirical sizes and empirical powers
		of the tests $T$, $T$, $T_{lc}$ and $T_{clx}$ in Case 1. These results are based on the 5\% significance level.}\label{tab1}
\end{table}

\begin{table}
\scriptsize
	\begin{tabular}{ccccccccccccc}
		\hline
		\multirow{3}{6.3em}{{\centering ($n_1,n_2,p$)\\$y_1>1$,$y_2>1$}}&\multicolumn{3}{c}{(25,35,40)} &\multicolumn{3}{c}{(50,70,80)} &\multicolumn{3}{c}{(100,140,160)} &\multicolumn{3}{c}{(200,280,320)}\\ \cline{2-13}
		&size&\multicolumn{2}{c}{power}
		&size&\multicolumn{2}{c}{power}&size&\multicolumn{2}{c}{power}&size&\multicolumn{2}{c}
		{power}\\\cline{2-13}
		&$a$=0&$a$=10&$a$=20&$a$=0 &$a$=10&$a$=20&$a$=0&$a$=10&$a$=20&$a$=0&$a$=10&$a$=20 \\ \hline
		$T$                     &0.052&0.961&1     & 0.046&0.953&  1   &0.046&0.950&   1  &0.048&0.948 & 1     \\ \hline
		$\tilde{T}$              &0.054&1    &1     & 0.050&1    &  1   &0.049&1    &   1  &0.048& 1    & 1     \\ \hline
		$T_{lc}$                &0.055&0.815&1     & 0.056&0.486&0.999 &0.050&0.209&0.959 &0.050&0.093 &0.579  \\ \hline
		$T_{clx}$                     &0.191&0.858&1     &0.125 &0.567&  1   &0.086&0.207&0.996 &0.070&0.095 &0.663  \\ \hline
		
		\multirow{3}{6.3em}{{\centering ($n_1,n_2,p$)\\$y_1>1$,$y_2<1$}}&\multicolumn{3}{c}{(25,35,30)} &\multicolumn{3}{c}{(50,70,60)} &\multicolumn{3}{c}{(100,140,120)} &\multicolumn{3}{c}{(200,280,240)}\\ \cline{2-13}
		&size&\multicolumn{2}{c}{power}
		&size&\multicolumn{2}{c}{power}&size&\multicolumn{2}{c}{power}&size&\multicolumn{2}{c}
		{power}\\\cline{2-13}
		&$a$=0&$a$=10&$a$=20&$a$=0 &$a$=10&$a$=20&$a$=0&$a$=10&$a$=20&$a$=0&$a$=10&$a$=20 \\ \hline
		$T$                     &0.056&0.701&0.999 &0.050 &0.610&0.998 &0.047&0.563&0.994 &0.053&0.517 &0.986  \\ \hline
		$\tilde{T}$              &0.057&0.999&1     & 0.052&1    &  1   &0.051&1    &   1  &0.050& 1    & 1     \\ \hline
		$T_{lc}$                &0.061&0.818&1     & 0.053&0.481&0.999 &0.050&0.205&0.956 &0.050&0.093 &0.573  \\ \hline
		$T_{clx}$                     &0.163&0.855&1     &0.112 &0.581&1     &0.087&0.228&0.997 &0.068&0.095 &0.679  \\ \hline
		\multirow{3}{6.3em}{{\centering ($n_1,n_2,p$)\\$y_1<1$,$y_2>1$}}&\multicolumn{3}{c}{(35,25,30)} &\multicolumn{3}{c}{(70,50,60)} &\multicolumn{3}{c}{(140,100,120)} &\multicolumn{3}{c}{(280,200,240)}\\ \cline{2-13}
		&size&\multicolumn{2}{c}{power}
		&size&\multicolumn{2}{c}{power}&size&\multicolumn{2}{c}{power}&size&\multicolumn{2}{c}
		{power}\\\cline{2-13}
		&$a$=0  &$a$=10 &$a$=20 &$a$=0  &$a$=10 &$a$=20 &$a$=0  &$a$=10 &$a$=20 &$a$=0  &$a$=10 &$a$=20  \\ \hline
		$T$                   &0.056	&0.108	&0.071	&0.052	&0.179	&0.272	&0.058	&0.225	&0.503	&0.055	&0.247	&0.654    \\ \hline
		$\tilde{T}$              &0.058	&0.977	&1	    &0.053	&0.993	&1	    &0.053	&0.997	&1	    &0.051	&0.998 	&1        \\ \hline
		$T_{lc}$                &0.058	&0.307	&0.995	&0.050	&0.137	&0.904	&0.052	&0.072	&0.488	&0.049	&0.051	&0.177    \\ \hline
		$T_{clx}$                     &0.169	&0.755	&1	    &0.117	&0.409	&0.998	&0.087	&0.156	&0.858	&0.066	&0.081	&0.303    \\ \hline
		\multirow{3}{6.3em}{{\centering ($n_1,n_2,p$)\\$y_1<1$,$y_2<1$}}&\multicolumn{3}{c}{(25,35,20)} &\multicolumn{3}{c}{(50,70,40)} &\multicolumn{3}{c}{(100,140,80)} &\multicolumn{3}{c}{(200,280,160)}\\ \cline{2-13}
		&size&\multicolumn{2}{c}{power}
		&size&\multicolumn{2}{c}{power}&size&\multicolumn{2}{c}{power}&size&\multicolumn{2}{c}
		{power}\\\cline{2-13}
		&$a$=0&$a$=10&$a$=20&$a$=0 &$a$=10&$a$=20&$a$=0&$a$=10&$a$=20&$a$=0&$a$=10&$a$=20 \\ \hline
		$T$                     &0.053&0.245&0.976& 0.055&0.117&0.724 &0.050&0.072&0.330 &0.051&0.055 &0.134    \\ \hline
		$\tilde{T}$              &0.061&0.999&1    & 0.056&1    &  1   &0.050&1    &   1  &0.047& 1    & 1     \\ \hline
		$T_{lc}$                &0.052&0.813&1    & 0.057&0.490&0.999 &0.052&0.204&0.957 &0.051&0.093 &0.568  \\ \hline
		$T_{clx}$                     &0.146&0.845&1    &0.106 &0.583&1     &0.078&0.237&0.996 &0.064&0.098 &0.706  \\ \hline
		
	\end{tabular}
	\caption{Empirical sizes and empirical powers
		of the tests $T$, $T$, $T_{lc}$ and $T_{clx}$ in Case 2. These results are based on the 5\% significance level.}\label{tab2}
\end{table}

\begin{table}[htbp]
\scriptsize
	\begin{tabular}{ccccccccccccc}
		\hline
		\multirow{3}{6.3em}{{\centering ($n_1,n_2,p$)\\$y_1>1$,$y_2>1$}}&\multicolumn{3}{c}{(25,35,40)} &\multicolumn{3}{c}{(50,70,80)} &\multicolumn{3}{c}{(100,140,160)} &\multicolumn{3}{c}{(200,280,320)}\\ \cline{2-13}
		&size&\multicolumn{2}{c}{power}
		&size&\multicolumn{2}{c}{power}&size&\multicolumn{2}{c}{power}&size&\multicolumn{2}{c}
		{power}\\\cline{2-13}
		&$a$=0&$a$=10&$a$=20&$a$=0 &$a$=10&$a$=20&$a$=0&$a$=10&$a$=20&$a$=0&$a$=10&$a$=20 \\ \hline
		$T$                     &0.051&0.960&1     & 0.046&0.954&  1   &0.046&0.949&   1  &0.048&0.949 & 1     \\ \hline
		$\tilde{T}$              &0.054&1    &1     & 0.050&1    &  1   &0.049&1    &   1  &0.048& 1    & 1     \\ \hline
		$T_{lc}$                &0.020&0.618&0.982 & 0.017&0.409&0.934 &0.016&0.219&0.746 &0.016&0.106 &0.455  \\ \hline
		$T_{clx}$                     &0.191&0.861&1     &0.125 &0.568&1     &0.086&0.210&0.996 &0.070&0.095 &0.662  \\ \hline
		
		\multirow{3}{6.3em}{{\centering ($n_1,n_2,p$)\\$y_1>1$,$y_2<1$}}&\multicolumn{3}{c}{(25,35,30)} &\multicolumn{3}{c}{(50,70,60)} &\multicolumn{3}{c}{(100,140,120)} &\multicolumn{3}{c}{(200,280,240)}\\ \cline{2-13}
		&size&\multicolumn{2}{c}{power}
		&size&\multicolumn{2}{c}{power}&size&\multicolumn{2}{c}{power}&size&\multicolumn{2}{c}
		{power}\\\cline{2-13}
		&$a$=0&$a$=10&$a$=20&$a$=0 &$a$=10&$a$=20&$a$=0&$a$=10&$a$=20&$a$=0&$a$=10&$a$=20 \\ \hline
		$T$                     &0.056&0.702&0.999 & 0.051&0.610&0.998 &0.046&0.564&0.994 &0.052&0.517 &0.986     \\ \hline
		$\tilde{T}$              &0.058&0.999&1     & 0.053&1    &  1   &0.052&   1 &   1  &0.049& 1    & 1     \\ \hline
		$T_{lc}$                &0.019&0.613&0.983 & 0.011&0.407&0.934 &0.016&0.216&0.755 &0.015&0.103 &0.458  \\ \hline
		$T_{clx}$                     &0.162&0.856&1     &0.112 &0.579&1     &0.087&0.229&0.997 &0.068&0.094 &0.679  \\ \hline
		\multirow{3}{6.3em}{{\centering ($n_1,n_2,p$)\\$y_1<1$,$y_2>1$}}&\multicolumn{3}{c}{(35,25,30)} &\multicolumn{3}{c}{(70,50,60)} &\multicolumn{3}{c}{(140,100,120)} &\multicolumn{3}{c}{(280,200,240)}\\ \cline{2-13}
		&size&\multicolumn{2}{c}{power}
		&size&\multicolumn{2}{c}{power}&size&\multicolumn{2}{c}{power}&size&\multicolumn{2}{c}
		{power}\\\cline{2-13}
		&$a$=0  &$a$=10 &$a$=20 &$a$=0  &$a$=10 &$a$=20 &$a$=0  &$a$=10 &$a$=20 &$a$=0  &$a$=10 &$a$=20  \\ \hline
		$T$                   &0.057	&0.106	&0.071	&0.052	&0.179	&0.274	&0.057	&0.223	&0.502	&0.054	&0.249	&0.648    \\ \hline
		$\tilde{T}$              &0.060	&0.978	&1	    &0.053	&0.992	&1	    &0.053	&0.997	&1	    &0.052	&0.998	&1        \\ \hline
		$T_{lc}$                &0.019	&0.313	&0.869	&0.016	&0.185	&0.667	&0.018	&0.098	&0.408	&0.017	&0.056	&0.220    \\ \hline
		$T_{clx}$                     &0.168	&0.757	&1	    &0.117	&0.407	&0.999	&0.088	&0.157	&0.855	&0.067	&0.081	&0.305    \\ \hline
		\multirow{3}{6.3em}{{\centering ($n_1,n_2,p$)\\$y_1<1$,$y_2<1$}}&\multicolumn{3}{c}{(25,35,20)} &\multicolumn{3}{c}{(50,70,40)} &\multicolumn{3}{c}{(100,140,80)} &\multicolumn{3}{c}{(200,280,160)}\\ \cline{2-13}
		&size&\multicolumn{2}{c}{power}
		&size&\multicolumn{2}{c}{power}&size&\multicolumn{2}{c}{power}&size&\multicolumn{2}{c}
		{power}\\\cline{2-13}
		&$a$=0&$a$=10&$a$=20&$a$=0 &$a$=10&$a$=20&$a$=0&$a$=10&$a$=20&$a$=0&$a$=10&$a$=20 \\ \hline
		$T$                     &0.053&0.250&0.975 & 0.054&0.119&0.725 &0.050&0.071&0.327 &0.052&0.054 &0.136    \\ \hline
		$\tilde{T}$              &0.060&0.999&1     & 0.056&   1 &  1   &0.050&1    &   1  &0.048& 1    & 1     \\ \hline
		$T_{lc}$                &0.019&0.629&0.983 & 0.016&0.403&0.930 &0.015&0.212&0.752 &0.015&0.113 &0.445  \\ \hline
		$T_{clx}$                     &0.147&0.847&1     &0.106 &0.581&1     &0.079&0.237&0.996 &0.063&0.098 &0.706  \\ \hline
		
	\end{tabular}
	\caption{Empirical sizes and empirical powers
		of the tests $T$, $T$, $T_{lc}$ and $T_{clx}$ in Case 3. These results are based on the 5\% significance level.}\label{tab3}
\end{table}

\begin{table}[htbp]
\scriptsize
	\begin{tabular}{ccccccccccccc}
		\hline
		\multirow{3}{6.3em}{{\centering ($n_1,n_2,p$)\\$y_1>1$,$y_2>1$}}&\multicolumn{3}{c}{(25,35,40)} &\multicolumn{3}{c}{(50,70,80)} &\multicolumn{3}{c}{(100,140,160)} &\multicolumn{3}{c}{(200,280,320)}\\ \cline{2-13}
		&size&\multicolumn{2}{c}{power}
		&size&\multicolumn{2}{c}{power}&size&\multicolumn{2}{c}{power}&size&\multicolumn{2}{c}
		{power}\\\cline{2-13}
		&$a$=0&$a$=10&$a$=20&$a$=0 &$a$=10&$a$=20&$a$=0&$a$=10&$a$=20&$a$=0&$a$=10&$a$=20 \\ \hline
		$T$                     &0.052&0.961&  1   & 0.046&0.954&  1   &0.045&0.950&   1  &0.047&0.949 & 1     \\ \hline
		$\tilde{T}$              &0.052&  1  &  1   & 0.049&1    &  1   &0.049& 1   &   1  &0.048& 1    & 1     \\ \hline
		$T_{lc}$                &0.083&0.505&0.909 & 0.081&0.348&0.793 &0.087&0.240&0.587 &0.082&0.179 &0.386  \\ \hline
		$T_{clx}$                     &0.062&0.207&0.778 &0.044 &0.076&0.492 &0.030&0.032&0.187 &0.022&0.020 &0.068  \\ \hline
		
		\multirow{3}{6.3em}{{\centering ($n_1,n_2,p$)\\$y_1>1$,$y_2<1$}}&\multicolumn{3}{c}{(25,35,30)} &\multicolumn{3}{c}{(50,70,60)} &\multicolumn{3}{c}{(100,140,120)} &\multicolumn{3}{c}{(200,280,240)}\\ \cline{2-13}
		&size&\multicolumn{2}{c}{power}
		&size&\multicolumn{2}{c}{power}&size&\multicolumn{2}{c}{power}&size&\multicolumn{2}{c}
		{power}\\\cline{2-13}
		&$a$=0&$a$=10&$a$=20&$a$=0 &$a$=10&$a$=20&$a$=0&$a$=10&$a$=20&$a$=0&$a$=10&$a$=20 \\ \hline
		$T$                     &0.055&0.701&0.999 & 0.050&0.611&0.998 &0.046&0.565& 0.994&0.053&0.518 &0.987    \\ \hline
		$\tilde{T}$              &0.057&0.999&1     & 0.053&  1  &  1   &0.052&1    &   1  &0.049& 1    & 1     \\ \hline
		$T_{lc}$                &0.089&0.517&0.918 & 0.083&0.359&0.786 &0.083&0.235&0.596 &0.081&0.172 &0.395  \\ \hline
		$T_{clx}$                     &0.068&0.230&0.818 &0.043 &0.085&0.539 &0.031&0.036&0.221 &0.021&0.023 &0.078  \\ \hline
		\multirow{3}{6.3em}{{\centering ($n_1,n_2,p$)\\$y_1<1$,$y_2>1$}}&\multicolumn{3}{c}{(35,25,30)} &\multicolumn{3}{c}{(70,50,60)} &\multicolumn{3}{c}{(140,100,120)} &\multicolumn{3}{c}{(280,200,240)}\\ \cline{2-13}
		&size&\multicolumn{2}{c}{power}
		&size&\multicolumn{2}{c}{power}&size&\multicolumn{2}{c}{power}&size&\multicolumn{2}{c}
		{power}\\\cline{2-13}
		&$a$=0  &$a$=10 &$a$=20 &$a$=0  &$a$=10 &$a$=20 &$a$=0  &$a$=10 &$a$=20 &$a$=0  &$a$=10 &$a$=20  \\ \hline
		$T $                   &0.057	&0.107	&0.069	&0.051	&0.178	&0.276	&0.057	&0.224	&0.506	&0.054	&0.248	&0.649    \\ \hline
		$\tilde{T}$              &0.060	&0.979	&1	    &0.054	&0.993	&1	    &0.053	&0.997	&1	    &0.051	&0.998	&1        \\ \hline
		$T_{lc}$                &0.084	&0.310	&0.736	&0.083	&0.229	&0.563	&0.081	&0.170	&0.375	&0.084	&0.125	&0.247    \\ \hline
		$T_{clx}$                     &0.065	&0.272	&0.777	&0.042	&0.145	&0.487	&0.029	&0.084	&0.244	&0.021	&0.048	&0.112    \\ \hline
		\multirow{3}{6.3em}{{\centering ($n_1,n_2,p$)\\$y_1<1$,$y_2<1$}}&\multicolumn{3}{c}{(25,35,20)} &\multicolumn{3}{c}{(50,70,40)} &\multicolumn{3}{c}{(100,140,80)} &\multicolumn{3}{c}{(200,280,160)}\\ \cline{2-13}
		&size&\multicolumn{2}{c}{power}
		&size&\multicolumn{2}{c}{power}&size&\multicolumn{2}{c}{power}&size&\multicolumn{2}{c}
		{power}\\\cline{2-13}
		&$a$=0&$a$=10&$a$=20&$a$=0 &$a$=10&$a$=20&$a$=0&$a$=10&$a$=20&$a$=0&$a$=10&$a$=20 \\ \hline
		$T$                     &0.052&0.245&0.975 & 0.055&0.119&0.724 &0.050&0.072&0.329 &0.051&0.055 &0.135   \\ \hline
		$\tilde{T}$              &0.061&0.999&1     & 0.056&1    &  1   &0.050&1    &   1  &0.048& 1    & 1     \\ \hline
		$T_{lc}$                &0.084&0.518&0.924 & 0.079&0.362&0.792 &0.083&0.249&0.595 &0.081&0.172 &0.398  \\ \hline
		$T_{clx}$                     &0.072&0.273&0.862 &0.043 &0.107&0.608 &0.032&0.044&0.261 &0.021&0.025 &0.092  \\ \hline
		
	\end{tabular}
	\caption{Empirical sizes and empirical powers
		of the tests $T$, $T$, $T_{lc}$ and $T_{clx}$ in Case 4. These results are based on the 5\% significance level.}\label{tab4}
\end{table}

\begin{table}\center{	\begin{tabular}{ccccccccc}
 \hline
		\multicolumn{3}{c}{(200,280)} &\multicolumn{3}{c}{(400,560)} &\multicolumn{3}{c}{(800,1120)} \\ \cline{1-9}
		p    &mean   &variance&p    &mean   &variance&p    &mean   &variance\\ \hline
		2    &0.0186 &0.1368  &4    &0.0099 &0.0441  &8    &0.0049 &0.0160  \\ \hline
		10   &0.0197 &0.0642  &20   &0.0081 &0.0260  &40   &0.0034 &0.0116  \\ \hline
		20   &0.0193 &0.0563  &40   &0.0107 &0.0252  &80   &0.0040 &0.0116  \\ \hline
		100  &0.0310 &0.0732  &200  &0.0121 &0.0338  &400  &0.0057 &0.0166  \\ \hline
		180  &0.0552 &0.1189  &360  &0.0223 &0.0555  &720  &0.0135 &0.0261  \\ \hline
		220  &0.0719 &0.1634  &440  &0.0332 &0.0736  &880  &0.0137 &0.0355  \\ \hline
		240  &0.0739 &0.1887  &480  &0.0395 &0.0870  &960  &0.0190 &0.0424  \\ \hline
		300  &0.1491 &0.3568  &600  &0.0695 &0.1598  &1200 &0.0373 &0.0771  \\ \hline
	\end{tabular}
	\caption{Numerical results for the estimator $\hat\Delta_1$ of standard normal distribution. Here the true value $\Delta_1=0$.}
\label{Del1}}
\end{table}

\begin{table}\center{	\begin{tabular}{ccccccccc}
\hline
		\multicolumn{3}{c}{(200,280)} &\multicolumn{3}{c}{(400,560)} &\multicolumn{3}{c}{(800,1120)} \\ \cline{1-9}
		p    &mean   &variance&p    &mean   &variance&p    &mean   &variance\\ \hline
		2    &-1.2021&0.0058  &4    &-1.2020&0.0032  &8    &-1.2012&0.0015  \\ \hline
		10   &-1.2009&0.0073  &20   &-1.2002&0.0036  &40   &-1.2018&0.0019  \\ \hline
		20   &-1.2005&0.0081  &40   &-1.2001&0.0039  &80   &-1.1994&0.0020  \\ \hline
		100  &-1.1971&0.0179  &200  &-1.1977&0.0091  &400  &-1.2001&0.0044  \\ \hline
		180  &-1.1858&0.0450  &360  &-1.1892&0.0213  &720  &-1.1957&0.0105  \\ \hline
		220  &-1.1704&0.0721  &440  &-1.1894&0.0314  &880  &-1.1938&0.0159  \\ \hline
		240  &-1.1650&0.0902  &480  &-1.1829&0.0417  &960  &-1.1892&0.0204  \\ \hline
		300  &-1.1170&0.2170  &600  &-1.1537&0.0974  &1200 &-1.1760&0.0445  \\ \hline
	\end{tabular}
	\caption{Numerical results for the estimator $\hat\Delta_1$ of uniform  distribution $(-\sqrt{3},\sqrt{3})$. Here the true value $\Delta_1=-1.2$.}
\label{Del2}}
\end{table}

\begin{table}\center{	\small	\centering
	\begin{tabular}{ccccc}
		\hline
		&Consumer Discretionary&Consumer Discretionary&Materials   &Materials \\
		
		\cline{2-5}
		($n_1$,$n_2$,$p$)   &{(61,62,71)}    &{(60,61,71)}     &{(60,58,22)}     &{(60,61,22)}  \\ \hline
		$T$            &0               &0                &0                &0                  \\ \hline
		$\tilde{T}$      &0               &$8.1307*10^{-8}$&0                &0
\\ \hline
		$T_{lc}$         &$9.4786*10^{-6}$&0                &0                &0           \\ \hline
		$T_{clx}$        &$4.7868*10^{-3}$ &$2.3442*10^{-4}$ &$1.9239*10^{-6}$ &$1.3389*10^{-3}$                   \\ \hline
	\end{tabular}
	\caption{$p$-values  of  the test of equality of the two covariance matrices of daily returns from the same sector in different  quarters.  }
	\label{realdata}}
\end{table}

\section{An example}
For illustration, we apply the proposed test statistics to the daily returns of
a selection of stocks issued by companies on Standard \& Poor's (S\&P) 500.
The original data are the closing prices or the bid/ask average of these stocks for the trading days during 2012 and 2013.  This dataset is derived from the Center for Research in Security Prices Daily Stock from Wharton Research Data Services. A common interest is to test whether  the covariance matrices of the  logarithmic daily returns for some stocks are the same over a period of time. Logarithmic daily returns are commonly used in finance. There are several theoretical and practical advantages of using logarithmic daily returns, including that we can assume that the sequences of logarithmic  daily returns are independent of each other, see \citep{Cont01E}.

According to the Global Industry Classification Standard (GICS), which is an industry taxonomy developed in 1999 by MSCI and S\&P for use by the global financial community, we select two  sectors, Consumer Discretionary and Materials, which include 71 and 22 stocks, respectively. {First, for Consumer Discretionary sector, we test whether the covariance matrices of the logarithmic daily returns of the second quarters (1 April-30 June) of 2012 and 2013 are the same. There were 63  and 64 trading days, respectively, in the  second quarters of 2012 and 2013. By using the logarithmic difference transformation on the original data, we obtain the   logarithmic daily returns dataset with sample sizes $n_1=61 $ and $n_2=62$. Simultaneously, we use the first quarters (1 January-31 March) data of 2012 and 2013 when applying to Materials sector. That makes the sample sizes be $n_1=60$ and $n_2=58$. The $p$ values obtained by applying the four test statistics  $T$, $\tilde{T}$,  $T_{lc}$ and $T_{clx}$ are shown in the first and third columns of Table \ref{realdata}.} Next, we use the same  procedure to test whether the covariance matrices of the logarithmic daily returns of the  first quarter (1 January-31 March) and the second quarter (1 April-30 June) of 2012  are the same. The results are  shown in the  second  and  fourth columns of Table \ref{realdata}. The results show that all the $p$-values are much smaller than 0.05. {Thus, there is strong evidence that the two covariance matrices  are different. Therefore, we  suggest  caution with the  assumption that the daily returns are identically distributed.}
%except $T_{clx}$ for the Materials sector. Thus, there is strong evidence that the two covariance matrices  are different. Therefore, we  suggest  caution with the  assumption that the daily returns are identically distributed.

\section{Conclusion  and discussion}
In this paper, we  propose two modified LRTs for  the equality of two high-dimensional covariance matrices and show the asymptotic distributions under the Moments Assumption and the null hypothesis. Furthermore, we show that the Moments Assumption and Dimensions Assumption are necessary for our results. More specifically, if the samples are Gaussian distributed, our modifications are optimal.  We also present the weakly consistent and asymptotic unbiased estimators of $\Delta_l$ for non-Gaussian distributions {under the null hypothesis}. According to the simulation results, the  performances of these  modified LRTs  are remarkable under conditions where they are applicable, but the theoretical results  for the alternative hypothesis are not considered in this paper because  of  the lack of random matrix theory. The  optimal modification of the LRT under  the alternative hypothesis test  will be presented in  the future.
\section*{Acknowledgement}
J. Hu was partially supported by  NSFC (No.11771073) and Science and Technology Development foundation of Jilin (No.20160520174JH). Z. D. Bai was partially supported by NSFC (No.11571067).

\appendix

\section{Appendix}

In this appendix, we provide the proofs of  Theorem \ref{th2} and Theorem \ref{th3}. By comparing  the definitions of $\mL$ and $\tilde \mL$, it is easy to verify that  the proof of Theorem \ref{th1} can be split into two parts,  one of which is  the same as the  proof of  Theorem \ref{th2} and another of which is an  analogous analysis.  Hence,  we present the proof {Theorem \ref{th2}} in this paper.

\subsection{Proof of Theorem \ref{th2}}

The main tool used to prove the theorems is the Cauchy integral formula and Theorem 1.6 in \cite{BaiH15C},  which established the central limit theorem of linear spectral statistics of random matrix $\bbB_n$  and is presented below for convenience.  Denote $\alpha_n=n_2/n_1$, $G_n(x)=p(\bbF^{\bbB_n}(x)-\bbF_{y_{1},y_{2}}(x))$ and $$\bbF_{y_{1},y_{2}}(x)=\frac{(\alpha_n+1)\sqrt{(x_r-x)(x-x_l)}}{2\pi y_1 x(1-x)}\delta_{x\in(x_l,x_r)}$$
$x_l,\:x_r=\frac{y_2(h\mp y_1)^2}{(y_1+y_2)^2}$
be the limit spectral distribution of $\bbB_n$ with parameters $\alpha_n$, $y_{1}$, $y_{2}$,
\begin{lemma}[Theorem 1.6 in \cite{BaiH15C}]\label{lemBH}
	
	In addition to the Moments Assumption and the Dimensions Assumption, we further assume that:
	
	(1) As $\min \{p, n_1, n_2\}\to\infty$, $y_1\to\gamma_1\in(0,1)\cup(1,\infty)$, $y_2\to\gamma_2\in(0,1)\cup(1,\infty)$, and $\alpha_n\to\alpha>0$.
	
	(2) Let $f_1,......f_k$ be the analytic functions on an open region containing the interval $[a_l,a_r]$, where $a_l=v^{-1}(1-\sqrt{\gamma_1})^2$, $a_r=1-\alpha v^{-1}(1-\sqrt{\gamma_2})^2$, and $v$ is defined as $v=(1+\frac{\gamma_1}{\gamma_2})(1-\sqrt{\frac{\gamma_1\gamma_2}{\gamma_1+\gamma_2}})^2$.
	
	Then, as min $(n_1,n_2,p)\rightarrow\infty$, the random vector $$(\int f_i dG_n(x)),\quad i=1,......,k,$$  converges weakly to a Gaussian vector $(G_{f_1},......G_{f_k})$ with mean functions
	\begin{equation*}
	\E G_{f_i}=
	\frac 1{4\pi i}\oint f_i(\frac{z}{\al+z})d\log(\frac{(1-\gamma_2)m_3^2(z)+2m_3(z)+1-\gamma_1}{(1+m_3(z))^2})
	\end{equation*}
	\begin{equation*}
	+\frac {\Delta_1}{2\pi i}\oint \gamma_1f_i(\frac {z}{\al+z})(1+m_3)^{-3}d m_3(z)
	\end{equation*}
	\begin{equation*}
	+\frac {\Delta_2}{4\pi i}\oint f_i(\frac {z}{\al+z})(1-\gamma_2m_3^2(z)(1+m_3(z))^{-2})d\log(1-\gamma_2m_3^2(z)(1+m_3(z))^{-2})
	\end{equation*}
	and covariance functions
	\begin{equation*}
	\textbf{Cov}(G_{f_i},G_{f_j})=-\frac{1}{2\pi^2}\oint\oint\frac{f_i(\frac{z_1}{\al+z_1})f_j(\frac{z_2}{\al+z_2})
		dm_3(z_1)dm_3(z_2)}{(m_3(z_1)-m_3(z_2))^2}
	\end{equation*}
	\begin{equation*}
	-\frac{\gamma_1\Delta_1+\gamma_2\Delta_2}{4\pi^2}\oint\oint\frac{f_i(\frac{z_1}{\al+z_1})f_j
		(\frac{z_2}{\al+z_2})dm_3(z_1)dm_3(z_2)}{(m_3
		(z_1)+1)^2(m_3(z_2)+1)^2},
	\end{equation*}
	where
\begin{align*}
	m_0(z)=&\frac{(1+\gamma_1)(1-z)-\alpha z(1-\gamma_2)+\sqrt{((1-\gamma_1)(1-z)+\alpha z(1-\gamma_2))^2-4\alpha z(1-z)}}{2z(1-z)(\gamma_1(1-z)+\alpha z\gamma_2)}-\frac 1{z},\\
	m_1(z)=&\frac{\alpha}{(\alpha+z)^2}m_0(\frac{z}{\alpha+z})-\frac1{\alpha+z} , \quad m_2(z)=-z^{-1}(1-\gamma_1)+\gamma_1m_1(z),\\
	m_{mp}^{\gamma_2}(z)=&\frac{1-\gamma_2-z+\sqrt{(z-1-\gamma_2)^2-4\gamma_2}}{2\gamma_2z},~m_3(z)=\gamma_2m_{mp}^{\gamma_2}(-m_2(z))+(m_2(z))^{-1}(1-\gamma_2).
\end{align*}
	All the above contour integrals can be evaluated on any contour enclosing the interval $[\frac{\alpha c_l}{1-c_l},\frac{\alpha c_r}{1-c_r}]$.
	
\end{lemma}
Note that Lemma \ref{lemBH} is proved based on the centralized sample covariance matrices, which are  constructed by not subtracting the sample mean vector from each sample vector. However, \cite{ZhengB15S} showed  that the only difference between the two types of  sample covariance matrices is normalization by $p/n_l
$ and $p/N_l$.
It is not difficult to verify that  Theorem \ref{th3}  satisfies Lemma \ref{lemBH} by choosing the kernel function to the logarithmic function.
Therefore, the main task of proving   Theorem \ref{th3} is to calculate the three  integrals, which will be shown step by step.
\begin{proof}[Proof of the limit part $\tilde \ell_n$]
	When calculating the integral \begin{align}\label{limp1}
	p\int\log{x}\frac{(\alpha_n+1)\sqrt{(x_r-x)(x-x_l)}}{2\pi y_1 x(1-x)}\delta_{x\in(x_l,x_r)}\mathrm{d}x
	\end{align}to achieve the limit part $\tilde \ell_n$, we choose a transformation $$x=\frac{y_2|y_1+h\xi|^2}{(y_1+y_2)^2}.$$ Clearly, when $x$ moves from $\frac{y_2(h- y_1)^2}{(y_1+y_2)^2}$ to $\frac{y_2(h+y_1)^2}{(y_1+y_2)^2}$ two times, $\xi$ shifts along a unit circle in the positive direction. Then, we obtain that the integral \eqref{limp1} is equal to
	\begin{align}\label{lneq1}
	p\frac{(y_1+y_2)h^2i}{4\pi}\oint  [\log\frac{y_2|y_1+h\xi|^2}{(y_1+y_2)^2}]
	\cdot\frac{(\xi^2-1)^2}{\xi^3|y_1+h\xi|^2|y_2-h\xi|^2}\mathrm{d}\xi
	\end{align}
	when $y_1>1$ and
	\begin{align}\label{lneq2}
	p\frac{(y_1+y_2)h^2i}{4\pi}\oint [ \log\frac{y_2|h+y_1\xi|^2}{(y_1+y_2)^2}]
	\cdot\frac{(\xi^2-1)^2}{\xi^3|y_1+h\xi|^2|y_2-h\xi|^2}\mathrm{d}\xi
	\end{align}
	when $y_1<1$. Two forms of the integral ensure that the logarithmic function returns a finite number of poles related to $y_1$ of the integrand. The pole related to $y_2$ of the integrand is $h/y_2$ when $y_2>1$. There is no differentce in the integral value before and after the transformation $\xi=1/\xi$, except that the residue point in the unit disc turns into $y_2/h$, which is the residue point under the assumption $y_2<1$. Thus, we assume $y_2>1$ without loss of generality. Then, we obtain that  if $y_1>1,$  \eqref{lneq1} can be rewritten as
	\bqn
	p\frac{(y_1+y_2)h^2i}{4\pi}\oint  [\log\frac{y_2(y_1+h\xi)}{(y_1+y_2)^2}]
	\times\frac{(\xi^2-1)^2}{\xi^3|y_1+h\xi|^2|y_2-h\xi|^2}\mathrm{d}\xi
	\eqn
	\bqn
	+p\frac{(y_1+y_2)h^2i}{4\pi}\oint  [\log\frac{y_2(y_1+\frac{h}{\xi})}{(y_1+y_2)^2}]
	\times\frac{(\xi^2-1)^2}{\xi^3|y_1+h\xi|^2|y_2-h\xi|^2}\mathrm{d}\xi
	\eqn
	which is equal to
	\bqn
	p\frac{(y_1+y_2)h^2i}{4\pi}\oint  [\log\frac{y_2(y_1+h\xi)^2}{(y_1+y_2)^2}]
	\times\frac{(\xi^2-1)^2}{\xi^3(y_1+h\xi)(y_1+\frac{h}{\xi})(y_2-h\xi)(y_2-\frac{h}{\xi})}\mathrm{d}\xi.
	\eqn  Similarly, if $y_1<1$, we have \eqref{lneq2} equals
	\bqn
	p\frac{(y_1+y_2)h^2i}{4\pi}\oint  [\log\frac{y_2(h+y_1\xi)^2}{(y_1+y_2)^2}]
	\times\frac{(\xi^2-1)^2}{\xi^3(y_1+h\xi)(y_1+\frac{h}{\xi})(y_2-h\xi)(y_2-\frac{h}{\xi})}\mathrm{d}\xi.
	\eqn  According to  Cauchy's residue theorem, we find three poles
	$$0,\quad -\frac{h}{y_1},\quad \frac{h}{y_2} $$ under the settings of $y_1>1, y_2>1.$ The corresponding residues are
	$$\frac{p(y_1+y_2)}{2y_1y_2}\log \frac{y_1^2y_2}{(y_1+y_2)^2},\quad\frac{p(1-y_1)}{2y_1}\log \frac{y_2(y_1-1)^2}{y_1^2},\quad\frac{p(1-y_2)}{2y_2}\log \frac{1}{y_2},$$ respectively. In the same way, under the assumptions $y_1<1, y_2>1,$ we obtain three poles
	$$0,\quad -\frac{y_1}{h},\quad \frac{h}{y_2}$$ and three residues
	$$\frac{p(y_1+y_2)}{2y_1y_2}\log \frac{h^2y_2}{(y_1+y_2)^2},\quad\frac{p(1-y_1)}{2y_1}\log \frac{h^2}{y_2(1-y_1)^2},\quad\frac{p(1-y_2)}{2y_2}\log \frac{h^2}{y_2}.$$
	Therefore, by combining the above results and basic calculations, we obtain the limit part  \eqref{limp1} is
	\begin{equation*}
	\log\left(\frac{ y_{2}h^{\frac{2h^2}{y_{1}y_{2}}}}{(y_{1}+y_{2})^{\frac{(y_{1}+y_{2})}{y_{1}y_{2}}}|1-y_{1}|^{\frac{|1-y_{1}|}{y_{1}}}}
	\right)
	-\log(\frac{{h^{\frac{2h^2}{y_{1}y_{2}}}}}{{y_1^\frac{1+y_2}{y_2}y_2^{\frac{1-y_1}{y_1}}}})\delta_{y_1>1}
	,\end{equation*}
	which completes the proof.
\end{proof}
\begin{proof}[Proof of the mean part $\tilde \mu_n$]
	Because  $m_3$ satisfies the equation $$z=-\frac{m_3(z)(m_3(z)+1-y_1)}{(1-y_2)(m_3(z)+\frac{1}{1-y_2})},$$ we make an integral conversion $z=(1+hr\xi)(1+\frac{h}{r\xi})/(1-y_2)^2$, where $r$ is a number greater than but close to 1. According to the discussion in the last section,  we assume $y_2>1$ without loss of generality. From the equation $$\frac{(1+hr\xi)(1+\frac{h}{r\xi})}{(1-y_2)^2}=-\frac{m_3(m_3+1-y_1)}{(1-y_2)m_3+1},$$ we obtain that $m_3=-(1+hr\xi)/(1-y_2)$ or $m_3=-(1+\frac{h}{r\xi})/(1-y_2)$. When $z$ runs in the positive direction around the contour enclosing the interval $[\frac{\alpha c_l}{1-c_l},\frac{\alpha c_r}{1-c_r}]$, $m_3$ runs in the opposite direction. Thus,  when $y_2>1$, we choose the outcome $m_3=-(1+\frac{h}{r\xi})/(1-y_2)$. Consequently, we have $$\frac{z}{\alpha+z}=\frac{y_2|1+hr\xi|^2}{|y_2+hr\xi|^2}.$$
	Therefore, for $y_1>1$, we get the mean part $\tilde \mu_n$ is equal to
	\begin{align}\label{mueq1}
	\lim_{r\downarrow1}\frac {1}{4\pi i}\oint_{|\xi|=1}(\log\frac{y_2(1+hr\xi)(1+\frac{h}{r\xi})}{(y_2+hr\xi)(y_2+\frac{h}{r\xi})})\times
	(\frac1{\xi-\frac1r}+\frac1{\xi+\frac1r}-\frac2{\xi+\frac h{y_2r}})\mathrm{d}\xi
	\end{align}
	\begin{align}\label{mueq2}
	+\lim_{r\downarrow1}\frac {\Delta_1}{2\pi i}\oint -y_1(\log\frac{y_2(1+hr\xi)(1+\frac{h}{r\xi})}{(y_2+hr\xi)(y_2+\frac{h}{r\xi})})
	\frac{(1-y_2)^2h}{y_2^3}\frac{\xi}{(\xi+\frac{h}{y_2r})^3}\mathrm{d}\xi
	\end{align}
	\begin{align}\label{mueq3}
	+\lim_{r\downarrow1}\frac {\Delta_2}{4\pi i}\oint(\log\frac{y_2(1+hr\xi)(1+\frac{h}{r\xi})}{(y_2+hr\xi)(y_2+\frac{h}{r\xi})})
	\frac{(y_2-1)(\xi^2-\frac{h^2}{y_2r^2})}{
		y_2(\xi+\frac{h}{y_2r})^2}[ \frac{2\xi}{(\xi^2-\frac{h^2}{y_2r^2})}-
	\frac{2}{\xi+\frac{h}{y_2r}}]\mathrm{d}\xi.
	\end{align}  If $y_1<1$,  we only  need to  change  the logarithmic term in the above expression   into the following form
	\bqn
	\log\frac{y_2(1+hr\xi)(1+\frac{h}{r\xi})}{(h+y_2r\xi)(h+\frac{y_2}{r\xi})},
	\eqn
	and the other terms remain the same.
	%and then we have $\tilde \mu_n$  equals to
	%\bqn
	%\lim_{r\downarrow1}\frac {1}{4\pi i}\oint_{|\xi|=1}(\log\frac{y_2(1+hr\xi)(1+\frac{h}{r\xi})}{(y_2+hr\xi)(y_2+\frac{h}{r\xi})})\times
	%(\frac1{\xi-\frac1r}+\frac1{\xi+\frac1r}-\frac2{\xi+\frac h{y_2r}})\mathrm{d}\xi
	%\eqn
	%\bqn
	%=\lim_{r\downarrow1}\frac {1}{4\pi i}\oint_{|\xi|=1}(\log\frac{y_2(1+hr\xi)}{(y_2+hr\xi)})\times
	%(\frac1{\xi-\frac1r}+\frac1{\xi+\frac1r}-\frac2{\xi+\frac h{y_2r}})\mathrm{d}\xi
	%\eqn
	%\bqn
	%+\lim_{r\downarrow1}\frac {1}{4\pi i}\oint_{|\xi|=1}(\log\frac{(1+\frac{h\xi}{r})}{(y_2+\frac{h\xi}{r})})\times
	%\frac{1}{\xi}\times(\frac{r}{r-\xi}+\frac{r}{r+\xi}-\frac{\frac{2y_2r}{h}}{\frac {y_2r}{h}+\xi})\mathrm{d}\xi.
	%\eqn
	For  $y_1>1$,  there are four poles of term \eqref{mueq1}
	$$\frac{1}{r},\quad -\frac{1}{r},\quad -\frac{h}{y_2r},\quad 0, $$
	and four residues
	$$\frac{1}{2}\log \frac{y_2(1+h)}{y_2+h},\quad \frac{1}{2}\log \frac{y_2(1-h)}{y_2-h},\quad \frac{1}{2}\log \frac{(y_1+y_2)^2}{y_2^2y_1^2},\quad 0.$$ Thus,  by  Cauchy's residue theorem,  we have
	\bqn
	\eqref{mueq1}
	=\frac{1}{2}\log\frac{(y_1+y_2)(y_1-1)}{y_1^2}.
	\eqn
	Analogously,
	\bqn
	\eqref{mueq2}
	=-\Delta_1\frac{h^2(y_1+y_2+y_1y_2)}{2y_1(y_1+y_2)^2},
	\eqn
	and
	\bqn
	\eqref{mueq3}
	=\Delta_2\frac{h^2y_2(2y_1^2-h^2)}{2y_1^2(y_1+y_2)^2}.
	\eqn
	For $y_1<1$,
	\bqn
	\lim_{r\downarrow1}\frac {1}{4\pi i}\oint_{|\xi|=1}(\log\frac{y_2(1+hr\xi)(1+\frac{h}{r\xi})}{(h+y_2r\xi)(h+\frac{y_2}{r\xi})})\times
	(\frac1{\xi-\frac1r}+\frac1{\xi+\frac1r}-\frac2{\xi+\frac h{y_2r}})\mathrm{d}\xi
	\eqn
	\bqn
	=\lim_{r\downarrow1}\frac {1}{4\pi i}\oint_{|\xi|=1}(\log\frac{y_2(1+\frac{h}{r\xi})}{(h+\frac{y_2}{r\xi})})\times
	(\frac1{\xi-\frac1r}+\frac1{\xi+\frac1r}-\frac2{\xi+\frac h{y_2r}})\mathrm{d}\xi
	\eqn
	\bqn
	+\lim_{r\downarrow1}\frac {1}{4\pi i}\oint_{|\xi|=1}(\log\frac{(\xi+hr)}{(h\xi+y_2r)})\times\frac{1}{\xi}
	\times(\frac{r}{r-\xi}+\frac{r}{r+\xi}-\frac{\frac{2y_2r}{h}}{\frac {y_2r}{h}+\xi})\mathrm{d}\xi
	\eqn
	%There are four poles
	%$$\frac{1}{r},\quad -\frac{1}{r},\quad -\frac{h}{y_2r},\quad 0,$$ and four residues
	%$$\frac{1}{2}\log\frac{y_2(1+h)}{h+y_2},\quad \frac{1}{2}\log \frac{y_2(1-h)}{h-y_2},\quad
	%\frac{1}{2}\log\frac{(y_1+y_2)^2}{y_2^2h^2},\quad 0.$$ Then we get
	%\bqn
	%\lim_{r\downarrow1}\frac {1}{4\pi i}\oint_{|\xi|=1}(\log\frac{y_2(1+hr\xi)(1+\frac{h}{r\xi})}{(h+y_2r\xi)(h+\frac{y_2}{r\xi})})\times
	%(\frac1{\xi-\frac1r}+\frac1{\xi+\frac1r}-\frac2{\xi+\frac h{y_2r}})\mathrm{d}\xi
	%\eqn
	\bqn
	=\frac{1}{2}\log\frac{(y_1+y_2)(1-y_1)}{h^2},
	\eqn
	\bqn
	\lim_{r\downarrow1}\frac {\Delta_1}{2\pi i}\oint -y_1(\log\frac{y_2(1+hr\xi)(1+\frac{h}{r\xi})}{(h+y_2r\xi)(h+\frac{y_2}{r\xi})})
	\frac{(1-y_2)^2h}{y_2^3}\frac{\xi}{(\xi+\frac{h}{y_2r})^3}\mathrm{d}\xi
	\eqn
	\bqn
	=\Delta_1\frac{y_1^2(y_1+2y_2)}{-2(y_1+y_2)^2},
	\eqn
	and
	\bqn
	\lim_{r\downarrow1}\frac {\Delta_2}{4\pi i}\oint(\log\frac{y_2(1+hr\xi)(1+\frac{h}{r\xi})}{(h+y_2r\xi)(h+\frac{y_2}{r\xi})})
	\frac{(y_2-1)(\xi^2-\frac{h^2}{y_2r^2})}{
		y_2(\xi+\frac{h}{y_2r})^2}[ \frac{2\xi}{(\xi^2-\frac{h^2}{y_2r^2})}-
	\frac{2}{\xi+\frac{h}{y_2r}}]\mathrm{d}\xi
	\eqn
	\bqn
	=\Delta_2\frac{y_1^2y_2}{2(y_1+y_2)^2}.
	\eqn
	Thus, by combining the above results,  the mean part $\tilde\mu_n$ is
	\bqn
	&\log\left[\frac{{(y_1+y_2)^\frac{1}{2}|1-y_1|^\frac{c_1}{2}}}{ h}\right]-\log(\frac{y_1}{h})\delta_{y_1>1}\\
	&-\frac{\Delta_1[y_1^3(y_1+2y_2)\delta_{y_1<1}+h^2(y_1+y_2+y_1y_2)\delta_{y_1>1}]}{2y_1(y_1+y_2)^2}
	+\frac{\Delta_2[y_1^4y_2\delta_{y_1<1}+h^2y_2(2y_1-h^2)\delta_{y_1>1}]}{2y_1^2(y_1+y_2)^2}.
	\eqn
	This completes the proof.
\end{proof}
\begin{proof}[Proof of the variance part]
	To calculate the variance part $\tilde \nu_n$
	\begin{equation}\label{nueq1}
	-\frac{1}{2\pi^2}\oint\oint\frac{f(\frac{z_1}{\al+z_1})f(\frac{z_2}{\al+z_2})
		dm_3(z_1)dm_3(z_2)}{(m_3(z_1)-m_3(z_2))^2}
	\end{equation}
	\begin{equation}\label{nueq2}
	-\frac{y_1\Delta_1+y_2\Delta_2}{4\pi^2}\oint\oint\frac{f(\frac{z_1}{\al+z_1})f
		(\frac{z_2}{\al+z_2})dm_3(z_1)dm_3(z_2)}{(m_3
		(z_1)+1)^2(m_3(z_2)+1)^2},
	\end{equation}
	we make an analogous integral conversion $$z_1=(1+hr_1\xi_1)(1+\frac{h}{r_1\xi_1})/(1-y_2)^2$$ and  $$z_2=(1+hr_2\xi_2)(1+\frac{h}{r_2\xi_2})/(1-y_2)^2.$$ Therefore, the relationship between $\xi_l$ and $m_3(z_l)$, $l=1$, $2$ is  $$m_3(z_1)=-\frac{1+\frac{h}{r_1\xi_1}}{(1-y_2)},~~~~m_3(z_2)=-\frac{1+\frac{h}{r_2\xi_2}}{(1-y_2)}.$$ Without loss of generality, we assume $r_1<r_2$.  When $y_1>1, y_2>1$,
	\bqn
	-\frac{1}{2\pi^2}\oint\oint\frac{f_i(\frac{z_1}{\al+z_1})f_j(\frac{z_2}{\al+z_2})
		dm_3(z_1)dm_3(z_2)}{(m_3(z_1)-m_3(z_2))^2}
	=2\lim_{r_2\downarrow1}\oint\frac{1}{2\pi i}
	(\log\frac{y_2(1+hr_2\xi_2)(1+\frac{h}{r_2\xi_2})}{(y_2+hr_2\xi_2)
		(y_2+\frac{h}{r_2\xi_2})})\times
	\eqn
	\bqn
	\{\lim_{r_1\downarrow1}\oint\frac{1}{2\pi i}
	(\log\frac{y_2(1+hr_1\xi_1)(1+\frac{h}{r_1\xi_1})}{(y_2+hr_1\xi_1)
		(y_2+\frac{h}{r_1\xi_1})})
	\frac{r_1}{(r_1\xi_1-r_2\xi_2)^2}\mathrm{d}\xi_1\}r_2\mathrm{d}\xi_2.
	\eqn
	\bqn
	=2\lim_{r_2\downarrow1}\oint\frac{1}{2\pi i}
	(\log\frac{y_2(1+hr_2\xi_2)(1+\frac{h}{r_2\xi_2})}{(y_2+hr_2\xi_2)
		(y_2+\frac{h}{r_2\xi_2})})\times
	\eqn
	\bqn
	\{\lim_{r_1\downarrow1}\oint\frac{1}{2\pi i}
	[(\log\frac{y_2(1+hr_1\xi_1)}{(y_2+hr_1\xi_1)
	})
	\frac{r_1}{(r_1\xi_1-r_2\xi_2)^2}+(\log\frac{(1+\frac{h}{r_1\xi_1})}{
		(y_2+\frac{h}{r_1\xi_1})})
	\frac{r_1}{r_2^2\xi_2^2(\frac{r_1}{r_2\xi_2}-\xi_1)^2}]\mathrm{d}\xi_1\}r_2\mathrm{d}\xi_2.
	\eqn
	There is only one pole $\frac{r_1}{r_2\xi_2}$ in the unit disc of the integration formula with respect to $\xi_1$. Thus, by    Cauchy's residue theorem, \eqref{nueq1} is equal to
	\bqn
	2\lim_{r_2\downarrow1}\oint\frac{1}{2\pi i}
	(\log\frac{y_2(1+hr_2\xi_2)(1+\frac{h}{r_2\xi_2})}{(y_2+hr_2\xi_2)
		(y_2+\frac{h}{r_2\xi_2})})\cdot\frac{1}{\xi_2}(\frac{h}{h+r_2\xi_2}-\frac{h}{h+y_2r_2\xi_2})
	\mathrm{d}\xi_2
	\eqn
	which has three poles $$0,\quad -\frac{h}{r_2},\quad -\frac{h}{r_2y_2}.$$
	Thus,
	we finally obtain that
	\bqn
	\eqref{nueq1} =\frac{2y_1^2}{(y_1+y_2)(y_1-1)}.
	\eqn
	Similarly, if $y_1<1,y_2>1$,
	%\bqn
	%-\frac{1}{2\pi^2}\oint\oint\frac{f_i(\frac{z_1}{\al+z_1})f_j(\frac{z_2}{\al+z_2})
	%	dm_3(z_1)dm_3(z_2)}{(m_3(z_1)-m_3(z_2))^2}
	%\eqn
	\bqn
	\eqref{nueq1} =\frac{2h^2}{(y_1+y_2)(1-y_1)}
	\eqn
	In addition, as when $y_1>1,y_2>1$,
	\begin{equation*}
	-\frac{y_1\Delta_1+y_2\Delta_2}{4\pi^2}\oint\oint\frac{f_i(\frac{z_1}{\al+z_1})f_j
		(\frac{z_2}{\al+z_2})dm_3(z_1)dm_3(z_2)}{(m_3
		(z_1)+1)^2(m_3(z_2)+1)^2}
	\end{equation*}
	\bqn
	=(y_1\Delta_1+y_2\Delta_2)\{\lim_{r_2\downarrow1}\oint\frac{1}{2\pi i}
	(\log\frac{y_2(1+hr_2\xi_2)(1+\frac{h}{r_2\xi_2})}{(y_2+hr_2\xi_2)
		(y_2+\frac{h}{r_2\xi_2})})
	\frac{(1-y_2)r_1h\mathrm{d}\xi_1}{(y_2r_1\xi_1+h)^2}\}\cdot
	\eqn
	\bqn
	\{\lim_{r_1\downarrow1}\oint\frac{1}{2\pi i}
	(\log\frac{y_2(1+hr_1\xi_1)(1+\frac{h}{r_1\xi_1})}{(y_2+hr_1\xi_1)
		(y_2+\frac{h}{r_1\xi_1})})
	\frac{(1-y_2)r_2h\mathrm{d}\xi_2}{(y_2r_2\xi_2+h)^2}\}
	\eqn
	\bqn
	=(y_1\Delta_1+y_2\Delta_2)\cdot\frac{h^4}{y_1^2(y_1+y_2)^2}
	\eqn
	and
	\begin{equation*}
	-\frac{y_1\Delta_1+y_2\Delta_2}{4\pi^2}\oint\oint\frac{f_i(\frac{z_1}{\al+z_1})f_j
		(\frac{z_2}{\al+z_2})dm_3(z_1)dm_3(z_2)}{(m_3
		(z_1)+1)^2(m_3(z_2)+1)^2}
	\end{equation*}
	\bqn
	=(y_1\Delta_1+y_2\Delta_2)\cdot\frac{y_1^2}{(y_1+y_2)^2}
	\eqn
	when $y_1<1,y_2>1.$ Thus, the variance part $\tilde\nu_n$ is equal to
	\begin{equation*}
	2\log\frac{h^2}{|1-y_1|(y_1+y_2)}+2\log(\frac{y_1^2}{h^2})\delta_{y_1>1}+\frac{(y_1\Delta_1+y_2\Delta_2)}
	{y_1^2(y_1+y_2)^2}[y_1^4\delta_{y_1<1}+h^4\delta_{y_1>1}].
	\end{equation*}
	This completes the proof.
\end{proof}
\subsection{Proof Theorem \ref{th3}}
We now  prove that  $\hat\Delta_1$ is a weakly consistent and asymptotically unbiased estimator of $\Delta_1$.
\begin{proof}
	From the definition of $\Delta_1$  in \eqref{delta1} and Chebyshev's inequality, we only need to prove the following two results: Under the same assumptions in Theorem \ref{th1} and  the null hypothesis,
	\begin{align}\label{th3eq1}
	\E \hat\Delta_1\to \Delta_1
	\end{align}
	and
	\begin{align}\label{th3eq2}
	\E (\hat\Delta_1- \Delta_1)^2\to 0.
	\end{align}
	We first consider \eqref{th3eq1}. Without loss of generality, we assume the mean of $\bbz_j^{(1)}$ is zero. {And  following the same  truncation steps in \citep{BaiS04C} we may truncate and re-normalize the random variables as follows}
\begin{align*}
|x_{ij}^{(1)}|\leq\eta_n\sqrt{n}, ~~\E x_{ij}^{(1)}=0, ~~\E (x_{ij}^{(1)})^2=1\quad\mbox{and} {\quad\E (x_{ij}^{(1)})^4=\Delta_1+3+O(n^{-1})},
\end{align*}
	where $\eta_n\to 0$ slowly. In addition, in the following we assume the sample covariance matrix without  the minus sample mean, because their difference is only a rank one matrix $\bar\bbz_1\bar\bbz_1'$ which will not affect the results. Form the  proof of  Theorem 1 in \citep{BaiY93L}, one can conclude that under the conditions of Theorem \ref{th1},  there exists a constant $M>0$ such that for any $k>0$
\begin{align*}
\mathbb{P}(\|(c_{11}\bbS^{\bbx}_{1j}+c_{12}\bbS_2^{\bbx})^{-1}\|\geq M)\leq n^{-k},
\end{align*}
where $\|\cdot\|$ means the spectral norm of a matrix. Thus, in the sequel, we can assume the smallest eigenvalue of $c_{11}\bbS^{\bbx}_{1j}+c_{12}\bbS_2^{\bbx}$ is bounded away from zero uniformly.
Then we have
%\begin{align*}
%\hat\Delta_1=(1-y)^2\frac{\sum_{j=1}^{N_1}[(\bbz^{(1)}_j-\bar{\bbz}^{(1)})'(c_{11}\bbS^{\bbz}_{1j}+c_{12}\bbS_2^{\bbz})^{-1}(\bbz^{(1)}_j-\bar{\bbz}^{(1)})-\frac{p}{1-y}]^2}{pN_1}
%-\frac{2}{1-y}\label{delta1}
%\end{align*}
\begin{align*}
&\E[(\bbz^{(1)}_j-\bar{\bbz}^{(1)})'(c_{11}\bbS^{\bbz}_{1j}+c_{12}\bbS_2^{\bbz})^{-1}(\bbz^{(1)}_j-\bar{\bbz}^{(1)})-\frac{p}{1-y}]^2\\
=&\E[(\bbz^{(1)}_j)'(c_{11}\bbS^{\bbz}_{1j}+c_{12}\bbS_2^{\bbz})^{-1}\bbz^{(1)}_j]^2+\frac{p^2}{(1-y)^2}\\
&-\frac{2p}{1-y}\E[(\bbz^{(1)}_j)'(c_{11}\bbS^{\bbz}_{1j}+c_{12}\bbS_2^{\bbz})^{-1}\bbz^{(1)}_j]+O(p).
\end{align*}
Here we use the fact that $\E[ (\bar{\bbz}^{(1)})'(c_{11}\bbS^{\bbz}_{1j}+c_{12}\bbS_2^{\bbz})^{-1}(\bar{\bbz}^{(1)})]=O(1)$ and $\E[ (\bar{\bbz}^{(1)})'(c_{11}\bbS^{\bbz}_{1j}+c_{12}\bbS_2^{\bbz})^{-1}({\bbz}_j^{(1)})]=O(1)$. These results can be found in \citep{PanZ11C}. According to the independence of $\bbz^{(1)}_j$ and $\bbS^{\bbz}_{1j}$ and under the null hypothesis,
we obtain that
\begin{align*}
\E[(\bbz^{(1)}_j)'(c_{11}\bbS^{\bbz}_{1j}+c_{12}\bbS_2^{\bbz})^{-1}\bbz^{(1)}_j]=\E[(\bbx^{(1)}_j)'(c_{11}\bbS^{\bbx}_{1j}+c_{12}\bbS_2^{\bbx})^{-1}\bbx^{(1)}_j]=\E tr [(c_{11}\bbS^{\bbx}_{1j}+c_{12}\bbS_2^{\bbx})^{-1}]
\end{align*}
and
\begin{gather*}
\E[(\bbz^{(1)}_j)'(c_{11}\bbS^{\bbz}_{1j}+c_{12}\bbS_2^{\bbz})^{-1}\bbz^{(1)}_j]^2
=(\Delta_1-3)\E tr[(c_{11}\bbS^{\bbx}_{1j}+c_{12}\bbS_2^{\bbx})^{-1}\circ (c_{11}\bbS^{\bbx}_{1j}+c_{12}\bbS_2^{\bbx})^{-1}]\\
+2\E tr[(c_{11}\bbS^{\bbx}_{1j}+c_{12}\bbS_2^{\bbx})^{-2}]+\{\E tr[(c_{11}\bbS^{\bbx}_{1j}+c_{12}\bbS_2^{\bbx})^{-1}]\}^2.
\end{gather*}
where $\circ$ is the Hadamard product. 
%Here we need some random matrix theories which can be found in \cite{BaiS10S} for more details. 
 Notice that
$c_{11}\bbS^{\bbx}_{1j}+c_{12}\bbS_2^{\bbx}$ is a sample covariance matrix with dimension $p$ and sample size $n_1+n_2-1$, whose  limit spectral distribution is the famous Marcenko-Pastur law with Stieltjes transform
\begin{align}\label{smp}
s_{mp}(z,\gamma)=\frac{1-\gamma-z+\sqrt{(z-1-\gamma)^2-4\gamma}}{2\gamma z},\quad z\in\mathbb{C}^+.
\end{align}
Here $ y:=p/(n_1+n_2-1)\to \gamma\in(0,\infty)$ and for  any distribution function $F$, its
Stieltjes transform  is defined  by
\begin{eqnarray*}
	s(z)=\int _{-\infty}^{+\infty}\frac{1}{x-z}dF(x) ,
\end{eqnarray*}
where $z\in \mathbb{C^{+}}$. Thus it is not difficult to check (see Section 3.3.2 in \cite{BaiS10S}) that
\begin{align}\label{deltaeq1}
\frac{1}{p}\E tr[(c_{11}\bbS^{\bbx}_{1j}+c_{12}\bbS_2^{\bbx})^{-1}]-\lim_{z\downarrow0}s_{mp}(z,y)\to0
\end{align}
and
\begin{align*}
\frac{1}{p}\E tr[(c_{11}\bbS^{\bbx}_{1j}+c_{12}\bbS_2^{\bbx})^{-2}]-\lim_{z\downarrow0}\frac{\partial s_{mp}(z,y)}{\partial z}\to0.
\end{align*}
From \eqref{smp}, L'Hospital's rule and the fact that $y<1$, we have
\begin{align}\label{deltaeq2}
\lim_{z\downarrow0}s_{mp}(z,y)-\frac{1}{1-y}=0
\end{align}
and
\begin{align*}{\lim_{z\downarrow0}\frac{\partial s_{mp}(z,y)}{\partial z}-\frac{1}{(1-y)^3}}=0.
\end{align*}
In addition, applying Lemma 4.3 of \cite{BaiH15C}, we obtain that
\begin{align*}
&\E tr[(c_{11}\bbS^{\bbx}_{1j}+c_{12}\bbS_2^{\bbx})^{-1}\circ (c_{11}\bbS^{\bbx}_{1j}+c_{12}\bbS_2^{\bbx})^{-1}]\\
=& \E\sum_{i=1}^p[(c_{11}\bbS^{\bbx}_{1j}+c_{12}\bbS_2^{\bbx})^{-1}(i,i)]^2=p^{-1}\{\E tr[(c_{11}\bbS^{\bbx}_{1j}+c_{12}\bbS_2^{\bbx})^{-1}]\}^2+o(p),
\end{align*}
which together with  \eqref{deltaeq1} and \eqref{deltaeq2} implies
\begin{align*}
\E tr[(c_{11}\bbS^{\bbx}_{1j}+c_{12}\bbS_2^{\bbx})^{-1}\circ (c_{11}\bbS^{\bbx}_{1j}+c_{12}\bbS_2^{\bbx})^{-1}]=\frac{p}{(1-y)^2}+o(p).
\end{align*}
Thus, combining the above results, we conclude that
\begin{align*}
\E[(\bbz^{(1)}_j-\bar{\bbz}^{(1)})'(c_{11}\bbS^{\bbz}_{1j}+c_{12}\bbS_2^{\bbz})^{-1}(\bbz^{(1)}_j-\bar{\bbz}^{(1)})-\frac{p}{1-y}]^2
=\frac{p{\Delta_1}}{(1-y)^2}+\frac{2p}{(1-y)^3} +o(p),
\end{align*}
which complete the proof of \eqref{th3eq1}.

For \eqref{th3eq2}, by the same argument as above, we have that $\E (\hat\Delta_1-\Delta_1)^2$ to
\begin{align*}
&\E (\hat\Delta_1-\Delta_1)^2\\
=&\frac{(1-y)^4\E\{\sum_{j=1}^{N_1}[{(\bbx^{(1)}_j)'}(c_{11}\bbS^{\bbx}_{1j}+c_{12}\bbS_2^{\bbx})^{-1}{\bbx^{(1)}_j}-tr(c_{11}\bbS^{\bbx}_{1j}+c_{12}\bbS_2^{\bbx})^{-1}]^2\}^2}{p^2N^2_1}\\&-{(\Delta_1+\frac{2}{1-y})^2}+o(1).
\end{align*}
Now, by (2.1) in \cite{BaiS04C},  we obtain that
\begin{align}\label{deltaeq3}
\frac{\sum_{j=1}^{N_1}\E[{(\bbx^{(1)}_j)'}(c_{11}\bbS^{\bbx}_{1j}+c_{12}\bbS_2^{\bbx})^{-1}
{\bbx^{(1)}_j}-tr(c_{11}\bbS^{\bbx}_{1j}+c_{12}\bbS_2^{\bbx})^{-1}]^4}{p^2N^2_1}=O(\eta_n^4{N_1^{-1}})\to0.
\end{align}
Next we use the fact that
\begin{align*}
(c_{11}\bbS^{\bbx}_{1j}+c_{12}\bbS_2^{\bbx})^{-1}
=&(c_{11}\bbS^{\bbx}_{1ji}+c_{12}\bbS_2^{\bbx})^{-1}\\&-\frac{{\frac{1}
{n_1+n_2-1}}(c_{11}{\bbS^
{\bbx}_{1ji}}+c_{12}\bbS_2^{\bbx})^{-1}{\bbx^{(1)}_i(\bbx^{(1)}_i)'}(c_{11}\bbS^{\bbx}_{1ji}+c_{12}\bbS_2^{\bbx})^
{-1}}{{1+\frac{1}{n_1+n_2-1}(\bbx^{(1)}_i)'(c_{11}\bbS^{\bbx}_{1ji}+c_{12}\bbS_2^{\bbx})^{-1}\bbx^{(1)}_i}}
\end{align*}
and
\begin{align*}
{(\bbx^{(1)}_i)'}(c_{11}\bbS^{\bbx}_{1j}+c_{12}\bbS_2^{\bbx})^{-1}{\bbx^{(1)}_i}=
\frac{{(\bbx^{(1)}_i)'}(c_{11}\bbS^{\bbx}_{1ji}+c_{12}\bbS_2^{\bbx})^{-1}
{\bbx^{(1)}_i}}{1+
{\frac{1}{n_1+n_2-1}
(\bbx^{(1)}_i)'}(c_{11}{\bbS^{\bbx}_{1ji}}+c_{12}\bbS_2^{\bbx})^{-1}{\bbx^{(1)}_i}},
\end{align*}
where  $\bbS^{\bbx}_{1ji}$ is the sample covariance matrix by removing the vector $\bbx_{j}^{(1)}$ and $\bbx_{i}^{(1)}$. Thus, we  get that
\begin{align*}
\sum_{i\neq j}^{N_1}\E&\{[{(\bbx^{(1)}_j)'}(c_{11}\bbS^{\bbx}_{1j}+c_{12}\bbS_2^{\bbx})^{-1}{\bbx^{(1)}_j}-tr(c_{11}\bbS^{\bbx}_{1j}+c_{12}\bbS_2^{\bbx})^{-1}]^2\\
&\times
[{(\bbx^{(1)}_i)'}(c_{11}{\bbS^{\bbx}_{1i}}+c_{12}\bbS_2^{\bbx})^{-1}{\bbx^{(1)}_i}-tr(c_{11}\bbS^{\bbx}_{1i}+c_{12}\bbS_2^{\bbx})^{-1}]^2\}
\end{align*}
\begin{gather*}
=N_1(N_1-1)\E\{[{(\bbx^{(1)}_1)'}(c_{11}\bbS^{\bbx}_{112}+c_{12}\bbS_2^{\bbx})^{-1}{\bbx^{(1)}_1}-tr(c_{11}\bbS^{\bbx}_{112}+c_{12}\bbS_2^{\bbx})^{-1}]^2\\
\times
[{(\bbx^{(1)}_2)'}(c_{11}\bbS^{\bbx}_{112}+c_{12}\bbS_2^{\bbx})^{-1}{\bbx^{(1)}_2}-tr(c_{11}\bbS^{\bbx}_{112}+c_{12}\bbS_2^{\bbx})^{-1}]^2\}+O(N_1^2)\\
=N_1(N_1-1)(\frac{p{\Delta_1}}{(1-y)^2}+\frac{2p}{(1-y)^3})^2+O(N_1^3),
\end{gather*}
which together with \eqref{deltaeq3} implies
\begin{align*}
&\frac{(1-y)^4\E\{\sum_{j=1}^{N_1}[{(\bbx^{(1)}_j)'}(c_{11}\bbS^{\bbx}_{1j}+c_{12}\bbS_2^{\bbx})^{-1}{\bbx^{(1)}_j}-tr(c_{11}\bbS^{\bbx}_{1j}+c_{12}\bbS_2^{\bbx})^{-1}]^2\}^2}{p^2N^2_1}\\
=&{(\Delta_1+\frac{2}{1-y})^2}+o(1).
\end{align*}
Then we complete  the proof of Theorem \ref{th3}.
\end{proof}

%
%\bibliographystyle{myabbrvnat}
%\bibliography{/Users/apple/Dropbox/huj/reference/paper,/Users/apple/Dropbox/huj/reference/book}
%%\bibliography{D:/Dropbox/huj/reference/paper,D:/Dropbox/huj/reference/book}
%%

\vskip .65cm
\noindent
Qiuyan Zhang, Jiang Hu and Zhidong Bai\\
School of Mathematics and Statistics\\
      Northeast Normal University,  China.
\vskip 2pt
\noindent
E-mail: (zhangqy919@nenu.edu.cn, huj156@nenu.edu.cn, baizd@nenu.edu.cn)

\end{document}